\documentclass[journal]{IEEEtran}
\usepackage{amssymb}
\usepackage{amsmath}
\usepackage{amsfonts}
\usepackage{graphicx}
\usepackage{cite}
\usepackage{blindtext, subfig}
\usepackage{color}
\usepackage{multicol}

\newtheorem{theorem}{Theorem}[section]
\newtheorem{assumption}{Assumption}[section]
\newtheorem{lemma}{Lemma}[section]

\newtheorem{remark}{Remark}[section]

\begin{document}
\graphicspath{ {images_20170214/} }  

\title{Global Finite-Time Attitude Consensus of Leader-Following Spacecraft Systems Based on Distributed Observers}

\author{Haichao~Gui and
        Anton~H.~J.~de~Ruiter
        \thanks{The authors are with the Department of Aerospace Engineering, Ryerson University, 350 Victoria Street, Toronto, ON, M5B 2K3, Canada.}
        }

\maketitle

\begin{abstract}
This paper addresses the leader-following attitude consensus problem for a group of spacecraft when at least one follower can access the leader's attitude and velocity relative to the inertial space.
A nonlinear distributed observer is designed to estimate the leader's states for each follower. The observer possesses one important and novel feature of keeping attitude and angular velocity estimation errors on second-order sliding modes, and thus provides finite-time convergent estimates for each follower. Further, quaternion-based hybrid homogeneous controllers recently developed for single spacecraft are extended and then applied, by establishing a separation principle with the proposed observer, to track the leader's attitude motion. As a result, global finite-time attitude consensus is achieved on the entire attitude manifold, with either full-state measurements or attitude-only measurements, as long as the network topology among the followers is undirected and connected. Numerical simulations are presented to demonstrate the performance of the proposed methods.
\end{abstract}



\section{Introduction}
Distributed attitude consensus of multiple cooperative spacecraft has drawn increasing attention due to its applications in formation flying, space-based interferometry, in-orbit assembly, etc. It can be classified as two types, namely, leaderless consensus that requires all spacecraft to reach an arbitrary yet probably \textit{a priori} unknown synchronized state \cite{Thunberg:14}, and leader-following consensus that requires each follower to track a prescribed group attitude trajectory provided by a real or virtual leader \cite{VanDyke:06,Abde:09,Dima:09}. This paper mainly focuses on the leader-following type.

The leader-following attitude consensus issue was first addressed by assuming that the leader's trajectory is available to all followers \cite{VanDyke:06,Abde:09,Dima:09,Mayhew:12,Ren:10}. In practice, a more common yet challenging case is that  only a portion of the followers can access the state of the leader. To deal with this problem, first-order sliding mode estimators were derived in \cite{Meng:10IJC,Zou:14} and \cite{Zou:16} to estimate the reference attitude and/or velocity in finite time. These designs can be traced back to the work of \cite{Cao:10} for single/double-integrator systems. Asymptotic distributed estimators were also proposed when the reference angular velocity is linearly parameterized \cite{Bai:08} or generated by a known, stable, linear system \cite{Cai:14,Cai:16}. The methods of \cite{Ren:07,Ren:10} and \cite{Du:11} involve no estimators but require the spacecraft to transmit their accelerations apart from their attitudes and velocities. Among the above methods, those of \cite{Du:11} and \cite{Zou:16}, further guarantee finite-time stability, which implies that the consensus behavior can be achieved in finite time instead of infinite time as for asymptotic or exponential stability. In addition, angular velocity measurements are not needed for the consensus algorithms of  \cite{Lawton:02,Abde:09,Ren:10,Zou:14,Zou:16} and \cite{Cai:16}.

Another important issue is the complex nonlinearity intrinsic in attitude control. More precisely, the attitude configuration, the set of $3\times 3$ rotation matrices SO(3), is a nonlinear manifold not diffeomorphic to any Euclidean space and precludes the existence of continuous, globally stabilizing, state-feedback laws on SO(3) \cite{Bhat:00,Chaturvedi:11}. In addition, the attitude kinematics and dynamics are both nonlinear. Due to these features, the attitude consensus laws extended from algorithms for linear systems ensure merely local or at most almost global stability \cite{Igarashi:09,Thunberg:14} while the methods  of \cite{Zou:14} and \cite{Zou:16} result in semi-global stability. In addition, some quaternion-based control schemes can cause the undesirable unwinding phenomenon due to neglecting that the unit-quaternion space is a double-covering of SO(3) \cite{Ren:07,Bai:08,Cai:14,Cai:16}. To overcome this problem, \cite{Mayhew:12} developed a hybrid feedback scheme with network-based hysteretic switching logics, resulting in global attitude  consensus and simultaneous robustness to measurement noise. This method, however, relies on the availability of the leader's state to all followers and, similarly to \cite{Abde:09} and \cite{Bai:08}, does not allow cycles in the communication graph. Otherwise, undesirable equilibria other than the consensus state can arise and fail the control objective.

This paper investigates the global attitude consensus of a leader-following spacecraft network in terms of the quaternion parameterization. The communication graph between followers is assumed to be an undirected connected graph and only a subset of the followers has access to the dynamic leader. In order to estimate the leader's states for each follower, a novel nonlinear distributed observer is designed such that finite-time convergence is guaranteed only if at least one follower connects to the leader. Following this, the hybrid homogeneous attitude controllers developed in \cite{Gui:16SCL} are extended and then applied together with the distributed observer to perform consensus control by establishing a separation principle \cite{Hong:00}. More precisely, the resultant consensus laws can restore the uniformly globally finite-time stable systems of \cite{Gui:16SCL}, in both the full-state measurement case and attitude-only measurement case, where the latter relies on a quaternion filter to inject the necessary damping instead of velocity feedback. As a result, the proposed control schemes avoid the unwinding problem and achieve global finite-time attitude consensus which, to the best of the our knowledge, has not been reported in existing cooperative attitude control literature. As another contribution, the proposed observer requires only the boundedness of the leader's angular velocity and its derivatives for finite-time convergence and hence possesses better robustness and allows more generic reference trajectories than the distributed observers in \cite{Cai:14,Cai:16} which are limited to stationary or periodic reference trajectories. In addition, it keeps the attitude and angular velocity estimation errors on second-order sliding modes, indicating higher accuracy during digital implementation than the distributed estimator derived in \cite{Cao:10} and its variants in \cite{Meng:10IJC,Zou:14} and \cite{Zou:16} that all attain first-order sliding modes.

The remainder of this paper is organized as follows. Basic notations, system equations and graph theory are reviewed in Section II. A nonlinear distributed observer with finite-time convergence is designed to estimate the leader's states in Section III. With the estimates from the observer, controller laws are derived in Section IV based on the results of \cite{Gui:16SCL} to attain global attitude consensus of the entire leader-following spacecraft system under two measurement scenarios. Section V presents numerical examples to illustrate the effectiveness of the proposed methods and Section VI draws the conclusions.

\section{Preliminaries}

Throughout this paper, denote by $ I_{n} $ the $ n\times n $ identity matrix, $\mathbf{1}_{n}=[1,...,1]^{T}$, and $\mathbb{I}_{n}= \{1,...,n\}$. For all $x\in{\mathbb R}$ and $\alpha \geq 0$, let $\textup{sgn}^{\alpha}(x)= \textup{sgn}(x)|x|^{\alpha}$ and $\textup{sat}_{\alpha}(x)= \textup{sgn}(x) \textup{min}\{ |x|^{\alpha}, 1\}$, where $\textup{sgn}(\cdot)$ is the standard sign function. Clearly, $\textup{sgn}^{\alpha}(x)$ is a continuous nonsmooth function if $0<\alpha<1$, while $\textup{sat}_{\alpha}(x)$ becomes the standard saturation function $\textup{sat}(x)$ if $\alpha = 1$. For all $x=[x_1, \cdots, x_n]^T \in{\mathbb R^{n}}$ and $\alpha \geq 0$, let $\textup{sgn}^{\alpha}(x)= [\textup{sgn}^{\alpha}(x_{1}), ..., \textup{sgn}^{\alpha}(x_{n})]$ and $\textup{sat}_{\alpha}(x)= [\textup{sat}_{\alpha}(x_{1}), ..., \textup{sat}_{\alpha}(x_{n})]$. Denote by $\|{\cdot}\|_{p}$ the $p$-norm of a vector respectively for $p=1,2,\infty$. For all $A\in\mathbb{R}^{m\times n}$, let $\bar{\sigma}(A)$ and $\underline{\sigma}(A)$ be its maximum and minimum singular values respectively. Note that $\bar{\sigma}(A)$ equals to its induced 2-norm $\|A\|_2= \textup{max}_{\{x\in\mathbb{R}^n: \|x\|_2=1\}} \|Ax\|_2$. Given $x\in {\mathbb R^{3}}$, $x^{\times}$ is the skew-symmetric matrix satisfying $x^{\times}y = x\times y$, $\forall y\in{\mathbb R^{3}}$, where $\times$ is the cross product on $\mathbb R^{3}$.

A quaternion $Q=[\eta,q^{T}]^{T} \in{\mathbb R^4}$ consists of a scalar part $\eta \in{\mathbb{R}}$ and a vector part $q\in \mathbb{R}^3$. Let $\textup{vec}(Q)$ give the vector part of $Q$, i.e., $\textup{vec}(Q)=q$. The quaternion multiplication is defined as
\begin{equation*}
 Q\circ Q^{\prime}=\left[
\begin{array}{c}
\eta{\eta^{\prime}} - q^{T}q^{\prime} \\
\eta q^{\prime} + \eta^{\prime}q + q\times q^{\prime}
\end{array}\right], Q^{\prime}=[\eta^{\prime}, q^{\prime T}]^T \in{\mathbb R^4}
\end{equation*}
which is associative and distributive but is not commutative. In addition, the conjugation of $Q$ is given by $Q^{*}= [\eta,-q^{T}]^T \in{\mathbb R^4}$. Note that $(Q\circ Q^{\prime})^{*}= (Q^{\prime})^{*}\circ Q^{*}$. A 3-D vector is treated as a quaternion with zero scalar part when operating with a quaternion. With the identity element $\textbf{1}=[1,0,0,0]^T$, the set of unit quaternions is defined as $\mathbb S^{3}= \{ Q\in{\mathbb R^4}: Q\circ Q^{*}=\textbf{1} \}$.

\subsection{System Equations}
Consider a system of \textit n rigid spacecraft (agents). Denote by $Q_{i}= [\eta_{i},q_{i}^{T}]^T \in{\mathbb S^3}$, $\forall i\in\mathbb{I}_n$, the attitude quaternion of the body-fixed frame of the \textit{i}th agent, $\mathcal F_{i}$, relative to the inertial frame $\mathcal F_{I}$. The equations of motion of the \textit{i}th agent are
\begin{equation} \label{eq_kine}
\dot{Q}_{i}= {1\over 2}Q_{i}\circ{\omega_{i}} = {1\over 2}E(Q_{i})\omega_{i}, \;
E(Q_{i})= \begin{bmatrix}
-q_{i}^T \\
q_{i}^{\times} + \eta_{i}I_3
\end{bmatrix}
\end{equation}
\begin{equation} \label{eq_dyn}
J_{i}\dot{\omega}_{i}= -\omega_{i}\times{J_{i}\omega_{i}} + u_{i},
\end{equation}
where $\omega_{i}\in{\mathbb R^3}$ and $J_{i}=J_{i}^{T}$  are the angular velocity and inertia tensor of agent \textit{i} expressed in $\mathcal F_i$. $u_i$ is the corresponding control torque. The rotation matrix from $\mathcal F_I$ to $\mathcal F_i$ can be computed from $Q_i$ by
\begin{equation} \label{eq_rotmat}
R(Q_i)= (\eta_i^2 - q_i^T q_i)I_3 - 2\eta_i q_{i}^{\times} + 2q_i q_i^T.
\end{equation}
Assume that the desired trajectory is generated by a leader spacecraft with body-fixed frame $\mathcal F_0$. Denote by $Q_0\in{\mathbb S^3}$ and $\omega_0\in{\mathbb R^3}$ the attitude quaternion and angular velocity of $\mathcal F_0$ relative to $\mathcal F_I$. In addition, $(Q_0,\omega_0)$ obeys the same kinematics as (\ref{eq_kine}). The attitude and angular velocity error of the \textit{i}th follower relative to the leader is defined as $Q_{i0}=Q_0^{*}\circ Q_i$ and $\omega_{i0}= \omega_i - R(Q_{i0})\omega_0$. Letting $\bar{\omega}_{i0}= R(Q_{i0})\omega_0$, the system equations in terms of $Q_{i0}$ and $\omega_{i0}$ are then written as
\begin{equation} \label{eq_errkine}
\dot{Q}_{i0}= {1\over 2}Q_{i0}\circ{\omega_{i0}} = {1\over 2}E(Q_{i0})\omega_{i0},
\end{equation}
\begin{equation} \label{eq_errdyn}
J_{i}\dot{\omega}_{i0}= \Xi(\omega_{i0},\bar{\omega}_{i0})\omega_{i0} -u_{fi} + u_{i},
\end{equation}
where $\Xi(\omega_{i0},\bar{\omega}_{i0})= (J_i\omega_{i0}+J_i\bar{\omega}_{i0})^{\times} - \bar{\omega}_{i0}^{\times}J_i - J_i\bar{\omega}_{i0}^{\times}$ is skew-symmetric and
\begin{equation} \label{eq_FF}
u_{fi}=J_{i}R(Q_{i0})\dot{\omega}_{0} + \bar{\omega}_{i0}^{\times}J_i\bar{\omega}_{i0},
\end{equation}
represents the torque to be compensated for perfect tracking of the desired trajectory. When every follower has access to the leader's trajectory, attitude consensus can be achieved by applying the controllers of \cite{Mayhew:11} and \cite{Gui:16SCL} to globally stabilize $(Q_{i0},\omega_{i0})=(\pm\textbf{1},0)$, $i\in{\mathbb I_n}$. These methods, however, cannot be applied if the leader's trajectory is available to only one or some of the followers.

As many studies on coordinated attitude control of formation flying spacecraft, the above dynamics models assume that all spacecraft share the same inertial frame $\mathcal{F}_I$. This is true in practice because the the inertial frame is usually set as the Earth-centered inertial frame for Earth spacecraft systems, and the heliocentric inertial frame for deep-space spacecraft systems.

\subsection{Communication Graph}

The information flow for \textit{n} followers is assumed to be bidirectional and can be described by a weighted undirected graph $\mathcal G \triangleq (\mathcal V, \mathcal E)$, where $\mathcal V = \{1,...,n\}$ is the node set and $\mathcal E \subseteq{\mathcal V \times \mathcal V}$ is the edge set. Since $\mathcal G$ is undirected, it follows that $(j,i)\in{\mathcal E} \Leftrightarrow (i,j)\in{\mathcal E}$, which means that there exists an edge between agents \textit{i} and \textit{j}. The adjacency matrix $A=[a_{ij}]_{n\times n}\in\mathbb{R}^{n\times n}$ is defined such that $a_{ij}=a_{ji}>0$ if $(j,i)\in{\mathcal E}$ while $a_{ij}=0$ otherwise. In addition, we set $a_{ii}=0$, $\forall i\in \mathcal V$. Denote by $L=[l_{ij}]_{n\times n} \in\mathbb{R}^{n\times n}$ the Laplacian matrix of $\mathcal G$ with $l_{ii}=\sum_{j=1}^{n}a_{ij}$ and $l_{ij}=-a_{ij}$ for $i\neq j$. Label the leader as node 0 and denote by $\bar{\mathcal G}$ the leader-following graph. Let $a_{i0}\geq 0$ be the connection weight between the leader and follower \textit{i} such that $a_{i0} = 1$ if follower \textit{i} connects to the leader and otherwise $a_{i0}=0$. Denote by $A_0= \textup{diag}\{a_{10},...,a_{n0}\}$. $\mathcal G$ is said to be connected if there is path between any two agents. Additionally, let $\bar {\mathcal G} \triangleq (\bar{\mathcal V}, \bar{\mathcal E})$ represent the leader-following graph, where $\bar{\mathcal V} = \{0\} \cup {\mathcal V}$ and $\bar{\mathcal E} \subseteq{\bar{\mathcal V} \times \bar{\mathcal V}}$. In order for the development of distributed attitude consensus schemes, the following assumptions and lemmas are introduced.
\begin{assumption} \label{assum_connect}
The communication among the followers is constant and bidirectional and the leader-following graph $\bar{\mathcal G}$ contains a spanning tree rooted at node 0, i.e., there is a path from the leader to any follower.
\end{assumption}
\begin{assumption} \label{assum_leader}                          
The leader's angular velocity $\omega_0(t)$ and its first two derivatives are continuous in time, and there exist constants $\gamma_{i}$, $i\in{\mathbb I_3}$, such that $\|\omega_0(t)\|_{\infty} \leq{\gamma_1}$, $\|\dot{\omega}_0(t)\|_{\infty} \leq{\gamma_2}$, and $\|\ddot{\omega}_0(t)\|_{\infty} \leq{\gamma_3}$.
\end{assumption}
\begin{lemma}  \label{lem_graph}                           
\cite{Hong:06} If Assumption \ref{assum_connect} holds, then $H\triangleq L+A_{0}$ is positive definite.
\end{lemma}

\begin{lemma}  \label{lem_STA}                           
\cite{Utkin:13} Consider the system $\dot{x}= - \rho_1\text{sgn}^{\frac{1}{2}}(x) + y$, $\dot{y} = -\rho_2 \text{sgn}(x) + f(t)$ with $|f(t)|\leq f_0$, where $f_0$, $\rho_1$, and $\rho_2$ are constants. If $\rho_1>0$ and $\rho_2>f_0$, then $x(t)$ and $y(t)$ converge to zero in a finite time $t_r \geq 0$, which can be estimated by the algorithm given in \cite{Utkin:13}.
\end{lemma}
\begin{lemma}  \label{lem_hardy}                           
\cite{Hardy:52} If $0<\alpha \leq 1$, then $(\sum_{i=1}^{n}|x_i|)^\alpha \leq \sum_{i=1}^{n}|x_i|^{\alpha} \leq n^{1-\alpha}(\sum_{i=1}^{n}|x_i|)^\alpha$ holds, $\forall x_i \in\mathbb R^n$ and $\forall i\in\mathbb I_n$.
\end{lemma}
The matrix $E(\cdot)$ defined in (\ref{eq_kine}) has the following useful properties and their proofs are given in Appendix \ref{app1}:
\begin{lemma} \label{lem_E1}                            
For any $Q,Q^{\prime}\in{\mathbb R^4}$, we have
$ \|E(Q)\|_2= \|Q\|_2$, ${Q^T}E(Q^{\prime})=-\textup{vec}^T(Q^*\circ{Q^\prime})$, and ${Q^T}E(Q)=0$.
\end{lemma}
\begin{lemma} \label{lem_E2}                            
For any $Q_i\in{\mathbb R^4}$, $b_{ij}\in{\mathbb R}$ and $b_{ij}=b_{ji}$, $i,j\in{\mathbb I_n}$ we have
$\sum_{i=1}^{n}\sum_{j=1}^{n} b_{ij}{Q_{j}^T}E(Q_i)=0.$  \end{lemma}

As shown in the next section, Lemmas \ref{lem_E1} and \ref{lem_E2} are useful for proving the stability of the proposed distributed observers.

\section{Distributed Finite-Time Observer Design}

In this section, a distributed observer is derived to obtain the leader's trajectory $(Q_{0},\omega_{0},\dot{\omega}_{0})$ for each agent in finite time when only a subset of the followers can access the leader's attitude and angular velocity relative to the inertial space.

Denote by $(P_i,v_i,z_i)\in\mathbb R^4 \times \mathbb R^3 \times \mathbb R^3$ the estimate of $(Q_{0},\omega_{0},\dot{\omega}_{0})$ by the \textit{i}th follower. Letting $(P_0,v_0,z_0)=(Q_{0},\omega_{0},\dot{\omega}_{0})$, a nonlinear distributed observer is then designed as
\begin{equation} \label{eq_obs1}
\dot{P}_i=\frac{1}{2} P_i\circ v_i -\lambda_1\textup{sgn}^{\beta_1}\left(\sum_{j=0}^{n}a_{ij}(P_i - P_j)\right ),
\end{equation}
\begin{equation} \label{eq_obs2}
\dot{v}_i= z_i -\lambda_2\textup{sgn}^{\beta_2}\left(\sum_{j=0}^{n}a_{ij}(v_i - v_j)\right ),
\end{equation}
\begin{equation} \label{eq_obs3}
\dot{z}_i= -\lambda_3\textup{sgn}\left(a_{i0}(z_i - w_i) + \sum_{j=1}^{n}a_{ij}(z_i - z_j) \right),
\end{equation}
where $(P_i(0),v_i(0),z_i(0))\in\mathbb R^4 \times \mathbb R^3 \times \mathbb R^3$, $i\in{\mathbb I_n}$, $\lambda_1,\lambda_2 >0$, $0<\beta_1,\beta_2<1$, and $\lambda_3>\gamma_3$. Since the leader's acceleration $\dot{\omega}_0$ is not available for any follower, it is recovered from the velocity $\omega_0$ via the second-order sliding mode differentiator
\begin{equation} \label{eq_obs4}
\left\{
\begin{array}{cll}
\dot{y}_i & = & -\mu_1 a_{i0} \text{sgn}^{\frac{1}{2}}(y_i - \omega_0) + w_i \\
\dot{w}_i & = & -\mu_2 a_{i0} \text{sgn}(y_i - \omega_0)
\end{array}\right.,
\end{equation}
where $y_i(0)= w_i(0) =0$, $\mu_1 >0$ and $\mu_2 \geq \gamma_3$. Note that the right-hand sides of $\dot{P}_i$, $\dot{v}_i$, and $\dot{y}_i$ are continuous while those of $\dot{z}_i$ and $\dot{w}_i$ are discontinuous and their solutions are understood in the Filippov sense \cite{Filippov:88}.

Next, we define the following estimation errors
\begin{equation} \label{eq_esterr}
\begin{array}{c}
\tilde{P}_i= P_i - P_0, \, \tilde{v}_i= v_i - v_0, \, \tilde{z}_i= z_i - z_0, \\
\tilde{y}_i= y_i - \omega_0, \, \tilde{w}_i= w_i - \dot{\omega}_0, \, i\in{\mathbb I_n}.
\end{array}
\end{equation}
In addition, denote by
\begin{multline} \label{eq_obsin}
x_{si}= \sum_{j=0}^{n}a_{ij}(x_i - x_j)= a_{i0}\tilde{x}_i + \sum_{j=1}^{n}a_{ij}(\tilde x_i - \tilde x_j), \\
x\in\{P,v,z\}, \, i\in{\mathbb I_n}.
\end{multline}
For convenience, introduce the aggregate variables
\begin{equation*}
\tilde{x}= [\tilde{x}_1^T,\cdots,\tilde{x}_n^T]^T, \,
x_{s}= [x_{s1}^T,\cdots,x_{sn}^T]^T, \, x\in\{P,v,z\},
\end{equation*}
which satisfy
\begin{equation} \label{eq_obsin_vec}
P_s = (H\otimes I_4)\tilde{P}, \, v_s= (H\otimes I_3)\tilde{v}, \, z_s= (H\otimes I_3)\tilde{z},
\end{equation}
where $\otimes$ is the Kronecker product. Differentiating the expressions in (\ref{eq_esterr}) and applying (\ref{eq_obs1})-(\ref{eq_obs4}), the equations of the estimation errors can then be written as
\begin{equation} \label{eq_errobs1}
\dot{\tilde{P}}_i=0.5(P_i\circ v_i - P_0\circ v_0) -\lambda_1\textup{sgn}^{\beta_1}(P_{si}), \, i\in{\mathbb I_n},
\end{equation}
\begin{equation} \label{eq_errobs2}
\dot{\tilde{v}}= \tilde{z} -\lambda_2\textup{sgn}^{\beta_2}(v_s),
\end{equation}
\begin{equation} \label{eq_errobs3}
\dot{\tilde{z}}= -\lambda_3\textup{sgn}(z_s - A_0\tilde{w}) - \textbf{1}_n\otimes{\dot{z}_0}.
\end{equation}
\begin{equation} \label{eq_errobs4}
\left\{
\begin{array}{cll}
\dot{\tilde{y}}_i & = & -\mu_1 a_{i0} \text{sgn}^{\frac{1}{2}}(\tilde{y}_i) + \tilde{w}_i \\
\dot{\tilde{w}}_i & = & -\mu_2 a_{i0} \text{sgn}(\tilde{y}_i) - \ddot{\omega}_0
\end{array}\right., \, i\in{\mathbb I_n},
\end{equation}
where $\tilde{w} = [\tilde{w}_1^T, \cdots, \tilde{w}_n^T]^T\in\mathbb{R}^{3n}$.

The following theorem shows that the proposed distributed observer ensures finite-time convergence of the estimation errors governed by (\ref{eq_errobs1})-(\ref{eq_errobs3}),
i.e., $(P_i(t),v_i(t),z_i(t))\rightarrow (Q_{0}(t),\omega_{0}(t),\dot{\omega}_{0}(t))$, $i\in{\mathbb I_n}$, in finite time.

\begin{theorem} \label{thm_obs}
Consider the distributed observer given by (\ref{eq_obs1})-(\ref{eq_obs4}) with  $\lambda_1,\lambda_2, \mu_1 >0$, $0<\beta_1,\beta_2<1$, and $\lambda_3, \mu_2>\gamma_3$. If Assumptions \ref{assum_connect} and \ref{assum_leader} hold, the state estimates are all uniformly bounded and $(P_i(t),v_i(t),z_i(t))=(Q_{0}(t),\\
\omega_{0}(t),\dot{\omega}_{0}(t))$, $i\in{\mathbb I_n}$, for all $t \geq T_p\geq 0$, where $T_p$ is given in (\ref{eq_Tp}).
\end{theorem}

\begin{IEEEproof}
See Appendix \ref{app2} for a detailed proof and an estimation of the convergence times for $P_i(t)$, $v_i(t)$, and $z_i(t)$, $i\in{\mathbb I_n}$, respectively. Notably, the proof also indicates that the convergence times can be made arbitrarily small by increasing $\lambda_{i}$, $i\in{\mathbb I_n}$.
\end{IEEEproof}

Note that a distributed asymptotic observer has recently been developed in \cite{Cai:14,Cai:16} for leader-following attitude consensus. This method is limited to the case that the leader's angular velocity is generated by a marginally stable linear system (i.e., $\dot{\omega}_0= S\omega_0$) and, particularly, the structure matrix $S$ must be precisely known by each follower. In contrast, the distributed observer derived here only requires the continuity and boundedness of the leader's trajectory (Assumption \ref{assum_leader}). Clearly, this condition is more generic and includes the marginally stable system $\dot{\omega}_0= S\omega_0$ as a special case, whether $S$ is known or not. In addition, the finite-time convergence property not only ensures high accuracy but also facilitates convenient verification of the separation principle between the proposed observer and many existing attitude controllers, as shown in the next section.

Since the dynamics of $P_i(t)$ and $v_i(t)$ given by (\ref{eq_obs1}) and (\ref{eq_obs2}) are continuous, the finite-time convergence property claimed in Theorem \ref{thm_obs} indicates that the identities $\tilde{P}_i(t)=\dot{\tilde{P}}_i(t)=0$ and $\tilde{v}_i(t)=\dot{\tilde{v}}_i(t)=0$, $i\in{\mathbb I_n}$, hold after a finite time. In other words, the proposed observer achieves second-order sliding modes for the attitude and angular velocity estimation errors, respectively. Distributed finite-time observers were constructed in \cite{Meng:10IJC,Zou:14} and \cite{Zou:16} to estimate the leader's attitude trajectory by extending the distributed sliding mode estimators developed for single/double-integrator systems \cite{Cao:10}. These methods, however, ensured merely first-order sliding modes for the attitude and angular velocity estimation errors because their derivatives involve discontinuous dynamics. Hence, our method is more desirable in the sense that higher-order sliding modes can produce better accuracy and less sensitivity to input noise during digital implementation, as shown in \cite{Levant:03}. Note also that the above distributed observer can be readily extended to the double-integrator systems and the Euler-Lagrange systems.

Similarly to \cite{Zou:16,Cai:14,Cai:16}, the attitude observer given by (\ref{eq_obs1}) evolves on $\mathbb{R}^4$ instead of $\mathbb{S}^3$. Therefore, it is possible that $P_i(t)$ is not a unit quaternion for $0\leq t \leq T_p$. In addition, the stability of (\ref{eq_obs1}) does not imply global stability on SO(3) and unwinding can occur in the proposed observer. Nonetheless, the influence of the possible unwinding can be mitigated in the sense that the observer convergence time $T_p$ can be made arbitrarily small by increasing $\lambda_i$, $i=1,2,3$. How to design global (asymptotic or finite-time) distributed attitude observers together with controllers on $\mathbb{S}^3$ or even directly on SO(3), like those for single spacecraft by \cite{Izadi:14,Bohn:16CDC,Bohn:16RNC}, remains an interesting open problem.

\section{Consensus Law Design}

Since the observer given by (\ref{eq_obs1})-(\ref{eq_obs3}) produces finite-time convergent estimates of the leader's trajectory for each follower, it can be combined with many existing attitude controllers developed for single spacecraft to reach group attitude consensus with the leader. To be effective, the salient fact that $P_i(t)$, $i\in\mathbb{I}_n$, is not necessarily a unit quaternion for $0\leq t < T_p$ , must be appropriately handled such that finite-time escape of the closed-loop trajectory does not occur during the observer transient. Next, we demonstrate how to deal with this issue by incorporating the observer with the global finite-time attitude controllers developed by \cite{Gui:16SCL} to solve the attitude consensus problem with full-state measurements and attitude-only measurements respectively.

\subsection{The Case of Full-State Measurements}

Assume that each follower spacecraft can measure its attitude and angular velocity relative to the inertial frame $\mathcal{F}_I$. Denote by $h_i\in\mathbb{H}\triangleq\{-1,1\}$ and $\mathbf{x}_{fi}=(Q_{i0},\omega_{i0},h_i)\in \mathcal{M}_1\triangleq\mathbb{S}^3\times \mathbb{R}^3\times \mathbb{H}$, $i\in\mathbb{I}_n$. Define $\overline{\textup{sgn}}:\mathbb{R}\to\mathbb{H}$ as an outer semicontinuous set-valued map, where $\overline{\textup{sgn}}(0)\in\mathbb{H}$ and $\overline{\textup{sgn}}(x)=\textup{sgn}(x)$ for $x\neq 0$. If follower $i$ has a direct connection to the leader, the following hybrid controller ensures uniform global finite-time stability of the equilibrium set $\mathcal{E}_i=\{\mathbf{x}_{fi}\in \mathcal{M}_1:Q_{i0}=h_i\textbf{1},\omega_{i0}=0\}$ \cite{Gui:16SCL}:
\begin{equation}\label{eq_GFTC1}
\left\{\begin{array}{lc}
\begin{split}
   u_i(\mathbf{x}_{fi})&=u_{fi}-k_{pi}\kappa_1(h_i Q_{i0},1-\alpha_{pi}) \\
     &\quad -k_{di}\textup{sat}_{\alpha_{di}}(\omega_{i0})
\end{split}   &\mathbf{x}_{fi}\in C_{fi}, \\
\mathbf{x}_{fi}^{+}=(Q_{i0},\omega_{i0}, \overline{\textup{sgn}}(\eta_{i0}))  &\mathbf{x}_{fi}\in D_{fi}.
\end{array}\right.
\end{equation}
where $k_{pi}, k_{di}>0$, $0<\alpha_{pi}<1$, $\alpha_{di}=2\alpha_{pi}/(1+\alpha_{pi})$, and $u_{fi}$ is given by (\ref{eq_FF}); $\mathbf{x}_{fi}^{+}$ denotes the state value immediately after a discontinuous jump and $\kappa_1(\cdot,\alpha)\in\mathbb{R}^3$ is given by
\begin{equation*}
\kappa_1(Q,\alpha)=\left \{ \begin{array}{cl}
                                   \frac{q}{\left(\sqrt{2(1 -\eta)}\right)^\alpha}, &\,\textup{if}\,\eta \neq 1,  \\
                                    0, & \,\textup{if}\,\eta=1,
                                 \end{array}\right.
\,Q\in\mathbb{S}^3.
\end{equation*}
The so-called \textit{flow set} $C_{fi}$ and \textit{jump set} $D_{fi}$ are defined as
\begin{equation}\label{eq_CD1}
\left\{\begin{split}
         C_{fi} &=\{\mathbf{x}_{fi}\in\mathcal{M}_1:h_i \eta_{i0} \geq -\delta\}, \\
           D_{fi}&= \{\mathbf{x}_{fi}\in\mathcal{M}_1:h_i \eta_{i0} \leq -\delta\}.
      \end{split}\right.
\end{equation}
where $0<\delta<1$. The above control law involves two control modes. More precisely, the continuous control torque $u_i(\mathbf{x}_{fi})$ is applied with $h_i$ unchanged (i.e., $\dot{h}_i=0$) when $\mathbf{x}_{fi}\in C_{fi}$ while if $\mathbf{x}_{fi}\in D_{fi}$, it jumps to $\mathbf{x}_{fi}^{+}$ immediately. Note that $(Q_{i0},\omega_{i0})$ remains continuous over the jump and only $h_i$ reverses its sign.

If follower \textit{i} does not have access to the leader, controller (\ref{eq_GFTC1}) cannot be applied since $Q_{i0}$, $\omega_{i0}$, and $u_{fi}$ are unavailable. By means of $(P_i(t),v_i(t),z_i(t))$ from (\ref{eq_obs1})-(\ref{eq_obs3}), estimates of $Q_{i0}$, $\omega_{i0}$, and $u_{fi}$ can be defined as

\begin{equation}\label{eq_errsig}
\hat{Q}_{i0}=[\hat{\eta}_{i0},\hat{q}_{i0}^T]^T= P_i^{*}\circ Q_i, \quad \hat{\omega}_{i0}= \omega_i - R(\hat{Q}_{i0})v_i,
\end{equation}
\begin{equation} \label{eq_estFF}
\hat{u}_{fi}=J_{i} R(\hat{Q}_{i0}) z_i + (R(\hat{Q}_{i0})v_i)^{\times}J_i R(\hat{Q}_{i0})v_i.
\end{equation}
It follows from Theorem \ref{thm_obs} that $\hat{Q}_{i0}(t)=Q_{i0}(t)$, $\hat{\omega}_{i0}(t)= \omega_{i0}(t)$, and $\hat{u}_{fi}(t)=u_{fi}(t)$ for $t \geq T_p$. Hence, the intuition motivates us to derive a control law from (\ref{eq_GFTC1}) by substituting $\hat{Q}_{i0}$, $\hat{\omega}_{i0}$, and $\hat{u}_{fi}$ respectively for $Q_{i0}$, $\omega_{i0}$, and $u_{fi}$. Such a design, however, can lead to unbounded control torques for $0\leq t \leq T_p$ and thus the instability of the closed-loop system. To see this, first note that the function $\kappa_1(\cdot,\alpha)$ for nonsmooth feedback injection is continuous and upper bounded (i.e., $\|\kappa_1(\cdot,\alpha)\|_2 \leq 1$) on $\mathbb{S}^3$ but not on $\mathbb{R}^4$ \cite{Gui:16SCL}. Since $\hat{Q}_{i0}(t)\in\mathbb{R}^4$ may not be a unit quaternion for $0\leq t \leq T_p$, it can occur, according to the definition of $\kappa_1(\cdot,\alpha)$, that $\kappa_1(\hat{Q}_{i0},\alpha)\to \infty$ as $\hat{\eta}_{i0}\to 1$ .

To overcome this problem, we define a function $\bar{\kappa}_1(\cdot,\alpha):\mathbb{R}^4 \to \mathbb{R}^3$ for $0\leq \alpha <1$ as

\begin{equation}\label{eq_kappa}
\bar{\kappa}_1(Q,\alpha)=\left \{ \begin{array}{cl}
                                   \frac{q}{\left(\sqrt{2\|Q\|_2(\|Q\|_2 -\eta)}\right)^\alpha}, &\,\textup{if}\,\eta \neq \|Q\|_2,  \\
                                    0, & \,\textup{if}\,\eta=\|Q\|_2.
                                 \end{array}\right.
\end{equation}
Clearly, $\bar{\kappa}_1(Q,\alpha)\to \kappa_1(Q,\alpha)$ as $\|Q\|_2 \to 1$. The following lemma further shows that $\bar{\kappa}_1(\cdot,\alpha)$ is continuous on $\mathbb{R}^4$, and its proof is given in Appendix \ref{app3}.

\begin{lemma} \label{lem_kappa}                            
The function $\bar{\kappa}_1(\cdot,\alpha)$ given by (\ref{eq_kappa}) is continuous on $\mathbb{R}^4$ and satisfies $\|\bar{\kappa}_1(Q,\alpha)\|_2 \leq \|Q\|_2^{1-\alpha}$, $\forall Q\in\mathbb{R}^4$.
\end{lemma}

Letting $\hat{\mathbf{x}}_{fi}= (\hat{Q}_{i0}, \hat{\omega}_{i0}, h_i) \in \hat{\mathcal{M}}_1 \triangleq \mathbb{R}^4 \times \mathbb{R}^3 \times \mathbb{H}$, the control law for follower \textit{i} is then designed as

\begin{equation}\label{eq_GFTC2}
\left\{\begin{array}{lc}
\begin{split}
   u_i(\hat{\mathbf{x}}_{fi})&=\hat{u}_{fi}- k_{pi}\bar{\kappa}_1(h_i \hat{Q}_{i0},1-\alpha_{pi}) \\
     & \quad -k_{di}\textup{sat}_{\alpha_{di}}(\hat{\omega}_{i0}),
\end{split} &\hat{\mathbf{x}}_{fi}\in \hat{C}_{fi}, \\
\hat{\mathbf{x}}_{fi}^{+}=(\hat{Q}_{i0},\hat{\omega}_{i0},
\overline{\textup{sgn}}(\hat{\eta}_{i0})), &\hat{\mathbf{x}}_{fi}\in \hat{D}_{fi}.
\end{array}\right.
\end{equation}
where the control parameters satisfy the same conditions as controller (\ref{eq_GFTC1}) and

\begin{equation}\label{eq_CD2}
\left\{\begin{split}
         \hat{C}_{fi} &=\{\hat{\mathbf{x}}_{fi}\in\hat{\mathcal{M}}_1:h_i \hat{\eta}_{i0} \geq -\delta\}, \\
           \hat{D}_{fi}&= \{\hat{\mathbf{x}}_{fi}\in\hat{\mathcal{M}}_1:h_i \hat{\eta}_{i0} \leq -\delta\}.
      \end{split}\right.
\end{equation}
Note that $u_i(\hat{\mathbf{x}}_{fi}(t))=u_i(\mathbf{x}_{fi}(t))$ for $t \geq T_p$. The following result then follows because the closed-loop system under controller (\ref{eq_GFTC2}) has no finite escape time.

\begin{theorem} \label{thm_GFTC2}
Consider the leader-following spacecraft system given by (\ref{eq_errkine}) and (\ref{eq_errdyn}). The hybrid control law given by (\ref{eq_GFTC2}) combined with the distributed observer given by (\ref{eq_obs1})-(\ref{eq_obs3}) ensures the uniform boundedness of $(Q_{i0}(t),\omega_{i0}(t))$ and globally stabilizes $(Q_{i0}(t),\omega_{i0}(t))$ to $(h_i\textbf{1},0)\in{\mathbb{S}^3\times \mathbb{R}^3}$, $i\in\mathbb{I}_n$, in finite time.
\end{theorem}

\begin{IEEEproof}
 See Appendix \ref{app4}.
\end{IEEEproof}

Compared to the quaternion-based, hybrid, asymptotic synchronization law derived in \cite{Mayhew:12}, the proposed consensus scheme not only ensures finite-time convergence but also removes the requirement of communicating the binary logic variable $h_i$ between neighboring agents, leading to a simpler switching logic. In addition, our method allows cyclic graph structure that can invalidate the method of \cite{Mayhew:12}.

As inferred by (\ref{eq_CD1}) and $0<\delta<1$, $h_i$ switches its sign in a hysteretic manner and only when the amount of sign mismatch between $h_i$ and $\eta_{i0}$ reaches the prespecified hysteresis width $\delta$. Large values of $\delta$ implies better robustness against measurement noise and less possibility of  undesirable chattering \cite{Mayhew:11}. Note that the hysteretic switching logic, if triggered, induces discontinuous command torques which is compatible with on/off thrusters. However, when implemented by actuators such as magnetic torquers, reaction wheels, and control moment gyros, the real output torque is a continuous approximation of the discontinuous jump and how accurate it can be depends on the actuator bandwidth. If the actuator responds ``fast enough'', the effect of hysteretic switching can be well approximated. A uniform bound on the maximum number of switching was established in \cite{Gui:16SCL} and shown to be proportional to the initial kinetic energy. Therefore, by initiating the control with small initial angular velocity errors can reduce the number of switching and thus avoid the occurrence of fast switching. More details are referred to \cite{Gui:16SCL}.

\subsection{The Case of Attitude-Only Measurements}

Next, consider the case that the \textit{i}th follower can obtain its attitude $Q_i$ but cannot measure its angular velocity $\omega_i$. Similarly to \cite{Tayebi:08}, the following quaternion filter is used to introduce damping for follower \textit{i}:
\begin{equation} \label{eq_QuatFilter}
\dot{\bar{Q}}_{i0}= {\frac{1}{2}}\bar{Q}_{i0}\circ{\bar{\Omega}_{i0}}, \quad \bar{Q}_{i0}(0)\in\mathbb{S}^3,
\end{equation}
where $\bar{\Omega}_{i0}\in \mathbb{R}^3$ is to be designed later. The quaternion error between $\hat{Q}_{i0}$ and $\bar{Q}_{i0}$ is given by $\tilde{Q}_{i0}= \bar{Q}_{i0}^{*}\circ{\hat{Q}_{i0}}$. By means of (\ref{eq_kine}), (\ref{eq_obs1}), (\ref{eq_errsig}), and (\ref{eq_QuatFilter}), the time derivative of $\tilde{Q}_{i0}$ is computed as
\begin{equation} \label{eq_DTildeQ}
\begin{split}
   \dot{\tilde{Q}}_{i0} &= \dot{\bar{Q}}_{i0}^*\circ \hat{Q}_{i0} + \bar{Q}_{i0}^*\circ \dot{\hat{Q}}_{i0} \\
     &={\frac{1}{2}} \tilde{Q}_{i0}\circ[\hat{\omega}_{i0} - R(\tilde{Q}_{i0}) \bar{\Omega}_{i0}] + {\frac{1}{2}}\Delta_i,
\end{split}
\end{equation}
where $R(\tilde{Q}_{i0})$ is computed from $\tilde{Q}_{i0}$ following (\ref{eq_rotmat}) and
\begin{equation*}
\begin{split}
\Delta_i &=(\|\tilde{Q}_{i0}\|_2^2 - 1)\bar{\Omega}_{i0}\circ\tilde{Q}_{i0} + \bar{Q}_{i0}^*\circ \\
&\quad\; [(\|\hat{Q}_{i0}\|_2^2 - 1)v_i\circ\hat{Q}_{i0} -2\lambda_1\textup{sgn}^{\beta_1}(P_{si}^*)\circ Q_i].
\end{split}
\end{equation*}
The presence of $\Delta_i$ in (\ref{eq_DTildeQ}) is due to the transient of the proposed observer. Noting $\|\tilde{Q}_{i0}(t)\|_2= \|\hat{Q}_{i0}(t)\|_2= \|P_i(t)\|_2$, it follows that $\Delta_i(t)=0$ for $t\geq T_p$.

Denote by $\hat{\mathcal{M}}_2\triangleq\mathbb{R}^4\times \hat{\mathcal{M}}_1\times \mathbb{H}$ and $\hat{\mathbf{x}}_{oi}=(\tilde{Q}_{i0},\hat{\mathbf{x}}_{fi},\tilde{h}_i)\in \hat{\mathcal{M}}_2$, where $\tilde{h}_i\in\mathbb{H}$ is the switching variable associated with $\tilde{Q}_{i0}$. The attitude control law for follower \textit{i} is designed as
\begin{equation}\label{eq_GFTC3}
\arraycolsep=2pt
\left\{\begin{array}{c}
        \left.\begin{array}{cll}
        u_i(\hat{\mathbf{x}}_{oi}) & = &\hat{u}_{fi}-k_{pi}\bar{\kappa}_1(h_i \hat{Q}_{i0},1-\alpha_{pi}) \\
        {}& {} & -k_{di}\bar{\kappa}_1(\tilde{h}_i \tilde{Q}_{i0},1-\alpha_{pi}) \\
        \bar{\Omega}_{i0} & =    & k_{qi} R^T(\tilde{Q}_{i0})\bar{\kappa}_1(\tilde{h}_i \tilde{Q}_{i0},1-\alpha_{qi})
        \end{array}
        \right\} \hat{\mathbf{x}}_{oi}\in \hat{C}_{oi}, \\
\hat{\mathbf{x}}_{oi}^{+}=(\tilde{Q}_{i0}, \hat{Q}_{i0},\hat{\omega}_{i0}, \overline{\textup{sgn}}(\hat{\eta}_{i0}), \overline{\textup{sgn}}(\tilde{\eta}_{i0})), \quad \hat{\mathbf{x}}_{oi}\in \hat{D}_{oi}.
\end{array}\right.
\end{equation}
where $k_{pi}, k_{di}, k_{qi}>0$, $0.5<\alpha_{qi}<1$, $\alpha_{pi}=2\alpha_{qi}- 1$. The flow set $\hat{C}_{oi}$ and jump set $\hat{D}_{oi}$ are given by

\begin{equation}\label{eq_CD3}
\left\{\begin{split}
         \hat{C}_{oi} &=\{\hat{\mathbf{x}}_{oi}\in\hat{\mathcal{M}}_2:h_i \hat{\eta}_{i0} \geq-\delta \: \textup{and} \: \tilde{h}_i\tilde{\eta}_{i0} \geq-\delta \}, \\
           \hat{D}_{oi}&= \{\hat{\mathbf{x}}_{oi}\in\hat{\mathcal{M}}_2:h_i \hat{\eta}_{i0} \leq-\delta \:\textup{or}\: \tilde{h}_i\tilde{\eta}_{i0} \leq-\delta  \}.
      \end{split}\right.
\end{equation}
When $\hat{\mathbf{x}}_{oi}\in \hat{C}_{oi}$, $h_i$ and $\tilde{h}_i$ remain unchanged. When $\hat{\mathbf{x}}_{oi}\in \hat{D}_{oi}$, the jump logic given in (\ref{eq_GFTC3}) only changes the sign of $h_i$ and/or $\tilde{h}_i$ while $(\tilde{Q}_{i0}(t), \hat{Q}_{i0}(t),\hat{\omega}_{i0}(t))$ remains continuous. The resultant closed-loop system by controller (\ref{eq_GFTC3}) is globally finite-time convergent, as stated in the following theorem.

\begin{theorem} \label{thm_GFTC3}
Consider the leader-following spacecraft system given by (\ref{eq_errkine}), (\ref{eq_errdyn}), and (\ref{eq_DTildeQ}). The hybrid control law given by (\ref{eq_GFTC3}) combined with the distributed observer given by (\ref{eq_obs1})-(\ref{eq_obs3}) ensures the uniform boundedness of $(\tilde{Q}_{i0}(t), Q_{i0}(t),\omega_{i0}(t))$ and globally stabilizes $(\tilde{Q}_{i0}(t),Q_{i0}(t),\omega_{i0}(t))$ to $(\tilde{h}_i\textbf{1}, h_i\textbf{1},0)\in{\mathbb{S}^3\times \mathbb{S}^3\times \mathbb{R}^3}$, $i\in\mathbb{I}_n$, in finite time.
\end{theorem}

\begin{IEEEproof}
For $t\geq T_p$, the fact of $(P_i(t),v_i(t),z_i(t))=(Q_{0}(t),\omega_{0}(t),\dot{\omega}_{0}(t))$, $i\in{\mathbb I_n}$, can be used to verify that the closed-loop equations (\ref{eq_errkine}), (\ref{eq_errdyn}), (\ref{eq_DTildeQ}) and (\ref{eq_GFTC3}) coincide with the uniformly globally finite-time stable hybrid system given in Theorem 3.3 of \cite{Gui:16SCL}. Hence, we only need to show that the $(Q_{i0}(t), \tilde{Q}_{i0}(t), \omega_{i0}(t))$ is bounded for $0\leq t<T_p$. Since $Q_{i0}(t)\in\mathbb{S}^3$ and $\|\tilde{Q}_{i0}(t)\|_2 = \|\hat{Q}_{i0}(t)\|_2= \|P_i(t)\|_2$ are both bounded, only the boundedness of $\omega_{i0}(t)$ needs to be verified. The details are analogous to the proof of Theorem \ref{thm_GFTC2} with obvious modifications and thus omitted.
\end{IEEEproof}

\begin{remark}
A distributed finite-time attitude consensus law has recently been developed in \cite{Zou:16} without angular velocity measurements. This method, however, is based on a non-global attitude representation and obtains merely semi-global stability. In addition, it relies on feedback input to cancel or at least dominate the entire nonlinear attitude dynamics, which can result in significant control expenditure. In contrast, our method avoids all these drawbacks and has a much simpler structure and thus better computational efficiency.
\end{remark}

\begin{remark}
When $\alpha_{pi}=\alpha_{di}=\alpha_{qi} =1$, the full-state feedback controller (\ref{eq_GFTC2}) and output-feedback controller (\ref{eq_GFTC3}) reduce to hybrid asymptotic controllers in \cite{Mayhew:11} for $t\geq T_p$, and thus lead to asymptotic convergence. Although the designs in this section mainly demonstrate the incorporation with the attitude controllers of \cite{Gui:16SCL}, one should keep in mind that many other (continuous) attitude controllers, e.g., \cite{Chen:09,deRuiter:13,Gui:17RNC} for single spacecraft can be integrated with the proposed distributed observer to perform group attitude consensus control so that different performance requirements are satisfied.
\end{remark}

\begin{remark}
The study of this paper assumes fixed communication topology and no time delays in the communication links and state measurements. Thunberg et al. \cite{Thunberg:14} addressed the leaderless attitude consensus with time-varying topologies in recent time but their study was built on kinematics level by treating angular velocities as inputs. How to design distributed attitude observers and controllers with switching topologies and time delays is a challenging topic for future research.
\end{remark}

\section{Numerical Simulations}

\begin{table*}
\renewcommand{\arraystretch}{1.5}
\centering
\caption{\label{tab1} Parameters for the distributed observer and controller}
\begin{tabular}{rc}
\hline\hline
{Algorithms} & {Parameters} \\
\hline
Observer~(\ref{eq_obs1})-(\ref{eq_obs4}) & $P_i(0)=Q_i(0)$, $v_i(0)=0$, $z_i(0)=[1,1,1]^T \, \text{rad}/\text{s}^2$, $i\in\mathbb{I}_n$, \\
{} & $\lambda_1=5$, $\lambda_2=1$, $\lambda_3=0.8$, $\beta_1=\beta_2=0.8$, $\mu_1 = 3$, $\mu_2=0.1$. \\
Controller~(\ref{eq_GFTC2}) & $k_{pi}=4$, $k_{di}=8$, $\alpha_{pi}=0.6$, $\alpha_{di}=0.75$, $\delta=0.2$, $h_i(0)=1$, $i\in\mathbb{I}_n$. \\
Controller~(\ref{eq_GFTC3}) & $k_{pi}=4$, $k_{di}=10$, $\alpha_{pi}=0.6$, $\alpha_{di}=0.75$, $\delta=0.2$, $h_i(0)=1$, \\
{} & $k_{qi}=3$, $\alpha_{qi}=0.8$, $\tilde{h}_i=1$,  $\bar{Q}_{i0}(0)=Q_i(0)$, $i\in\mathbb{I}_n$. \\
\hline\hline
\end{tabular}
\end{table*}

\begin{figure}
\begin{center}
\includegraphics[height=4cm]{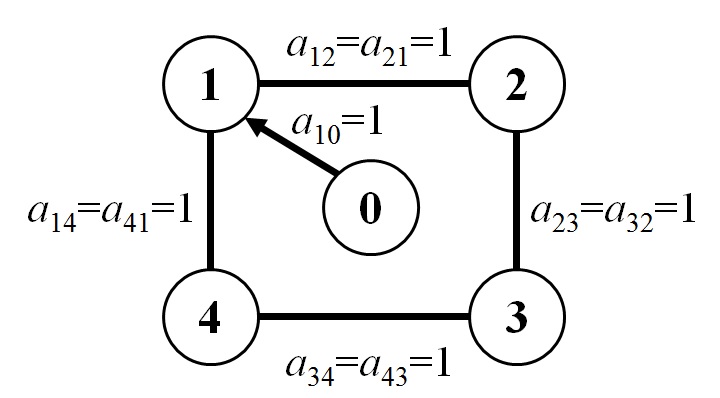}
\caption{Communication graph of the leader-following system}
\label{fig1}
\end{center}
\end{figure}

In the following, the performance of the proposed distributed control methods are demonstrated through a simulated leader-following spacecraft system with four identical rigid spacecraft as followers. The leader's attitude trajectory is described by $Q_0(0)=\textbf{1}$ and $\omega_0(t)= 0.01[\sin(\Omega_0 t),\cos(\Omega_0 t),\sin(\Omega_0 t)]^T$ rad/s, where $\Omega_0= 0.01$ rad/s. The inertia of the follower spacecraft is given by $J_i=\textup{diag}(10,8,12) \, \textup{kg}\cdot\textup{m}^2$, $i\in\mathbb{I}_n$. Three numerical examples are presented, namely, 1) simulations without any uncertainty, 2) simulations with communication and measurement delays as well as external disturbances, and 3) simulations with time-varying topologies. The simulations are conducted in MATLAB/Simulink using the default integrator ode45 with a maximum step size of 0.001 s.

\subsection{Simulations without Any Uncertainty}

In this section, simulations are conducted with zero disturbance and instantaneous communication between neighboring spacecraft. Assume that the communication graph $\bar{G}$ is fixed and all nonzero edge weights are given in Fig. \ref{fig1}. The measurement and communication frequencies of each spacecraft are both set to 100 Hz. The initial attitudes and angular velocities of the four followers are given by

\begin{equation*}
\begin{matrix}
\begin{bmatrix}
  Q_1^T(0) \\
  Q_2^T(0) \\
  Q_3^T(0) \\
  Q_4^T(0)
\end{bmatrix} & = & \begin{bmatrix}
                      0 & 1 & 0 & 0 \\
                      0 & 0 & -1 & 0 \\
                      0.6164 & 0.5 & -0.6 & 0.1 \\
                      -0.8426 & -0.2 & 0.3 & 0.4
                    \end{bmatrix}
\end{matrix}
\end{equation*}

\begin{equation*}
\begin{matrix}
\begin{bmatrix}
  \omega_1^T(0) \\
  \omega_2^T(0) \\
  \omega_3^T(0) \\
  \omega_4^T(0)
\end{bmatrix} & = & \begin{bmatrix}
                      0.2 & 0.2 & 0.2 \\
                      -0.1 & -0.1 & -0.1 \\
                      0.4 & 0.4 & 0.4 \\
                      -0.3 & -0.3 & 0.3
                    \end{bmatrix}
\end{matrix} \, \textup{rad/s}
\end{equation*}

\begin{figure*}
\begin{center}
\subfloat[]{\includegraphics[height=4cm] {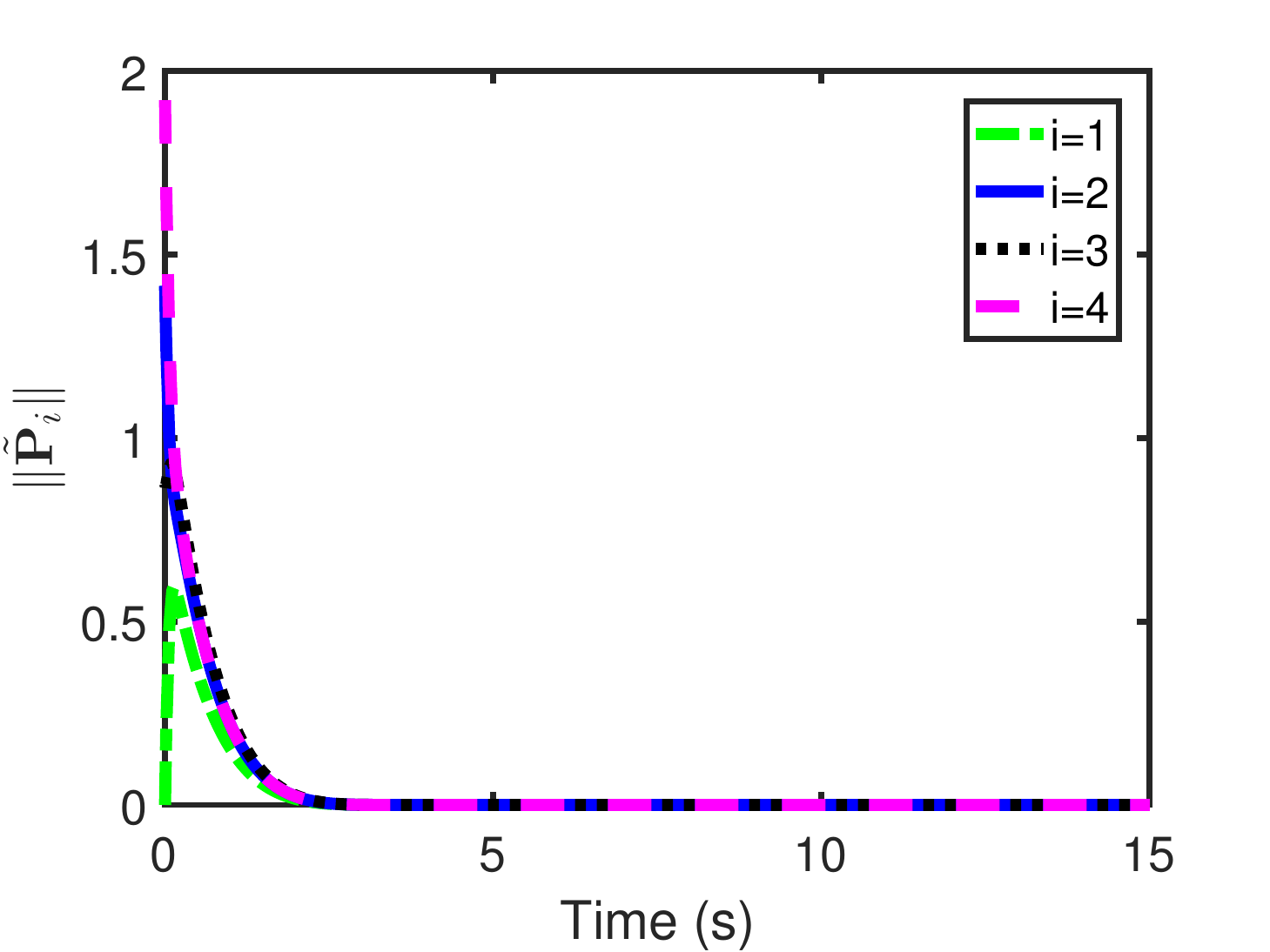} \label{fig2a}}
\subfloat[]{\includegraphics[height=4cm] {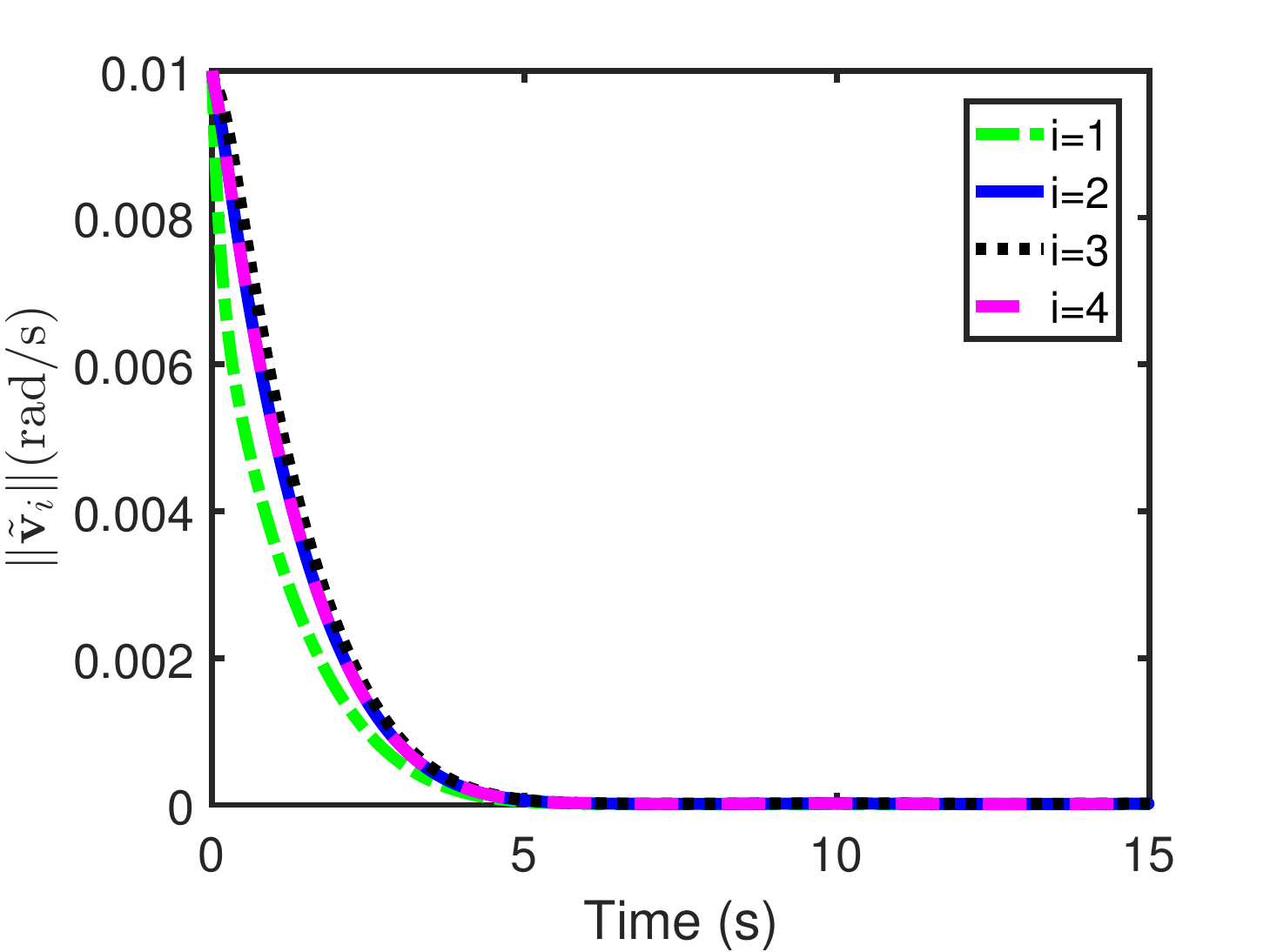} \label{fig2b}}
\subfloat[]{\includegraphics[height=4cm] {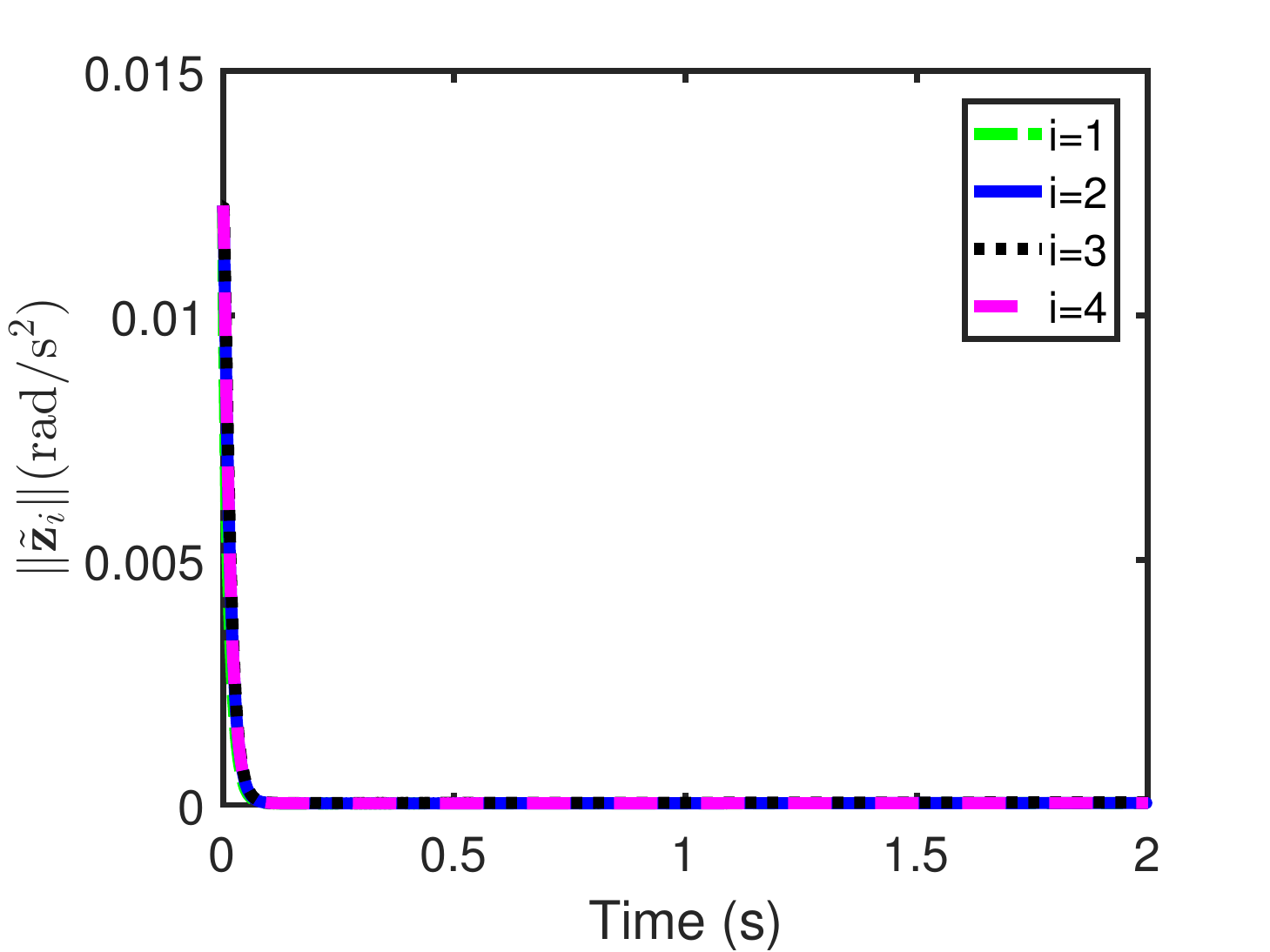} \label{fig2c}}
\caption{Response of the distributed observer without any uncertainty: (a) $\|\tilde{P}_i\|_2$, (b) $\|\tilde{v}_i\|_2$, and (c) $\|\tilde{z}_i\|_2$}  \label{fig2}
\end{center}
\end{figure*}

\begin{figure*}
\begin{center}
\subfloat[]{\includegraphics[height=4cm] {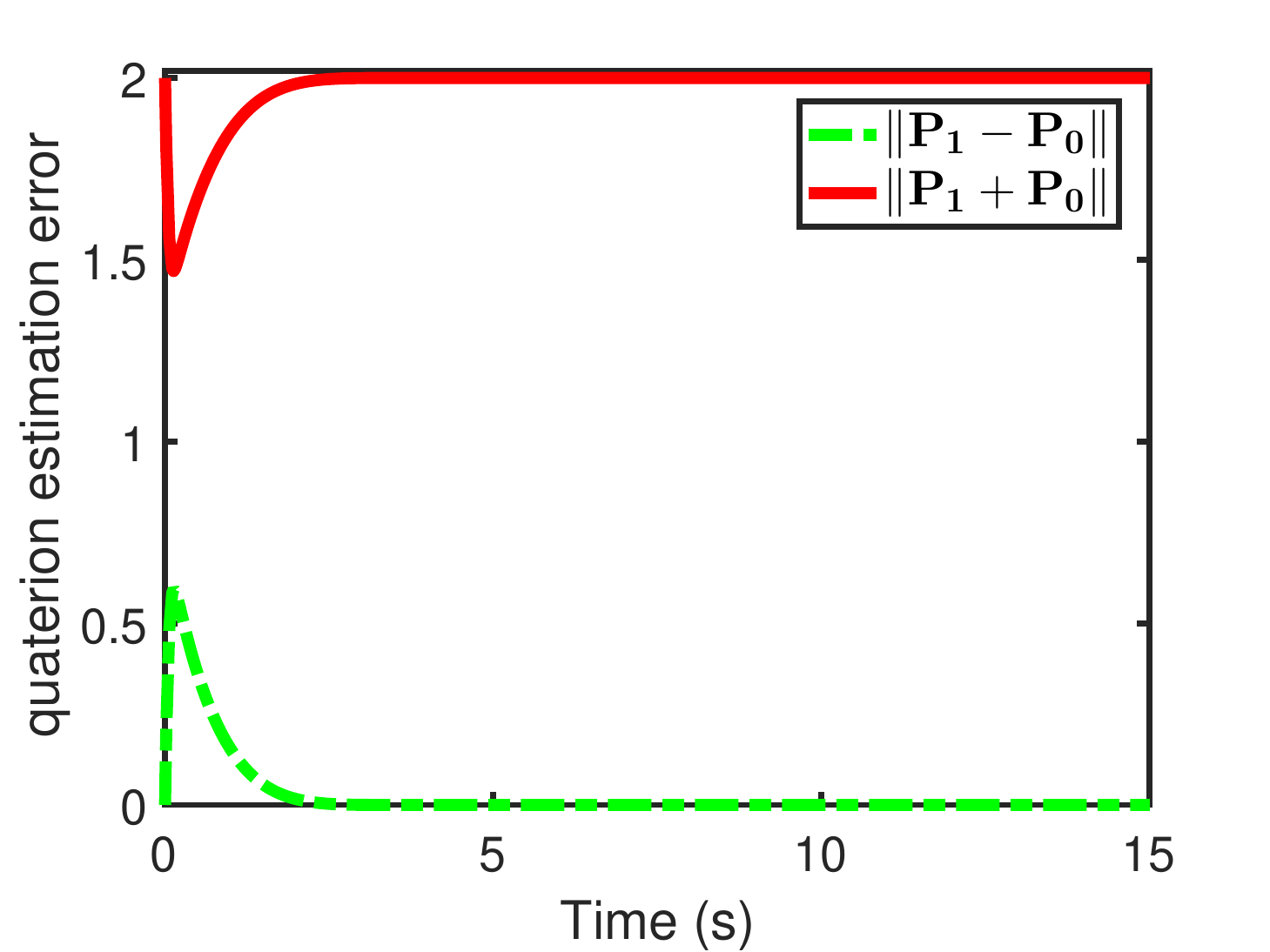}}
\subfloat[]{\includegraphics[height=4cm] {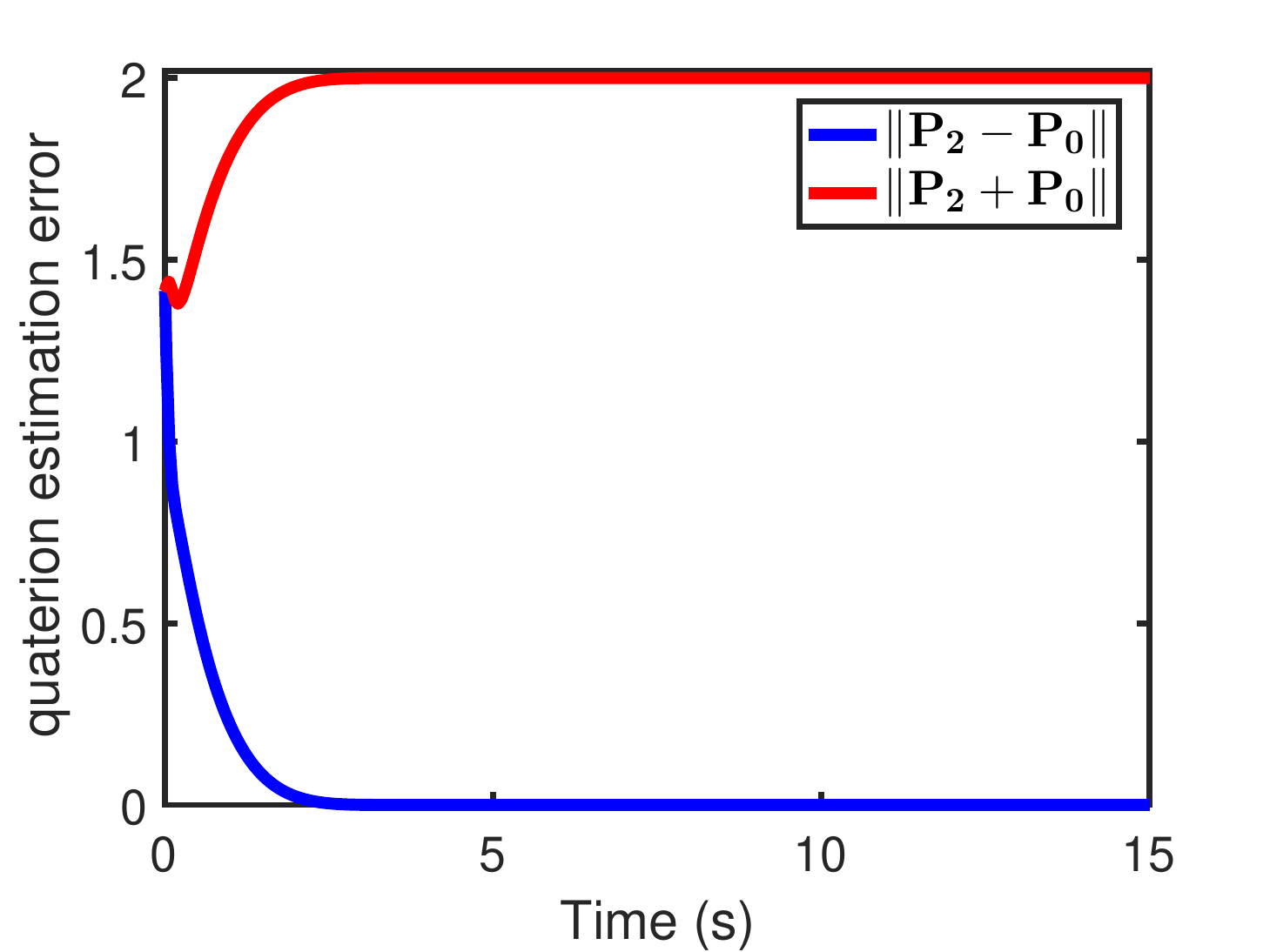}} \\
\subfloat[]{\includegraphics[height=4cm] {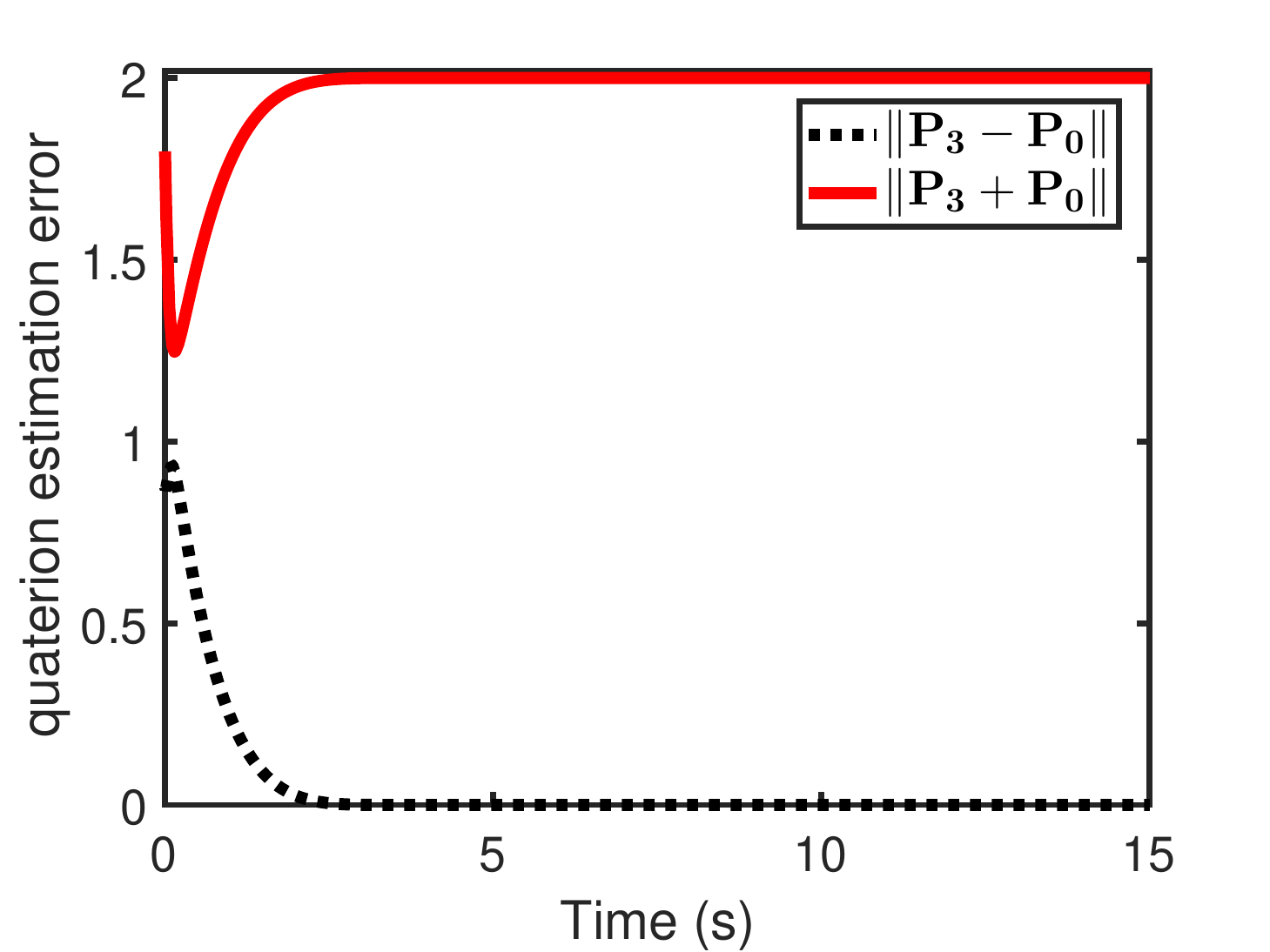}}
\subfloat[]{\includegraphics[height=4cm] {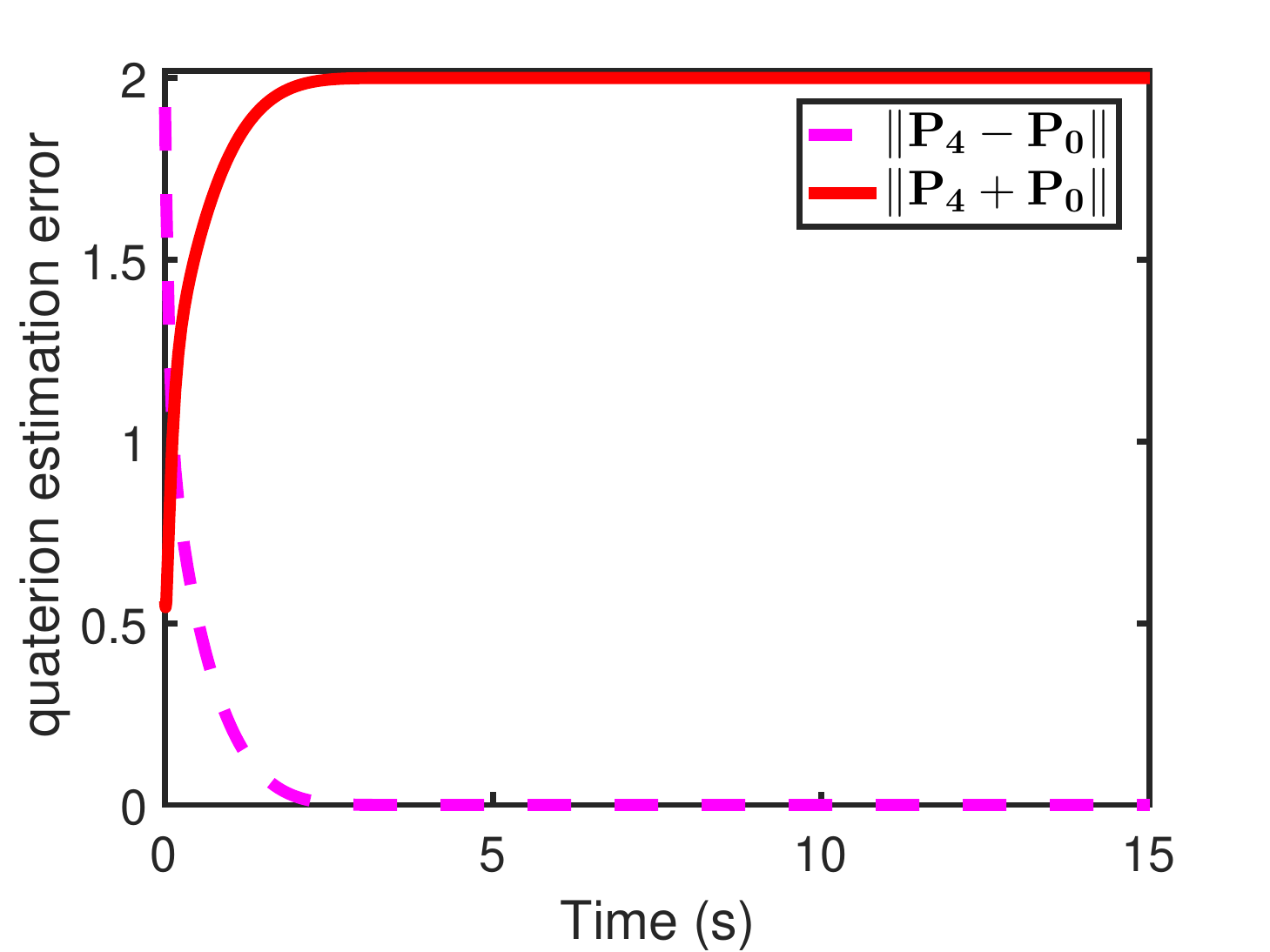} }
\caption{The leader's attitude estimation error $\|P_i - P_0\|$ and $\|P_i + P_0\|$ , $i\in\mathbb{I}_4$}  \label{fig:unwinding}
\end{center}
\end{figure*}

\begin{figure*}
\begin{center}
\subfloat[]{\includegraphics[height=4cm] {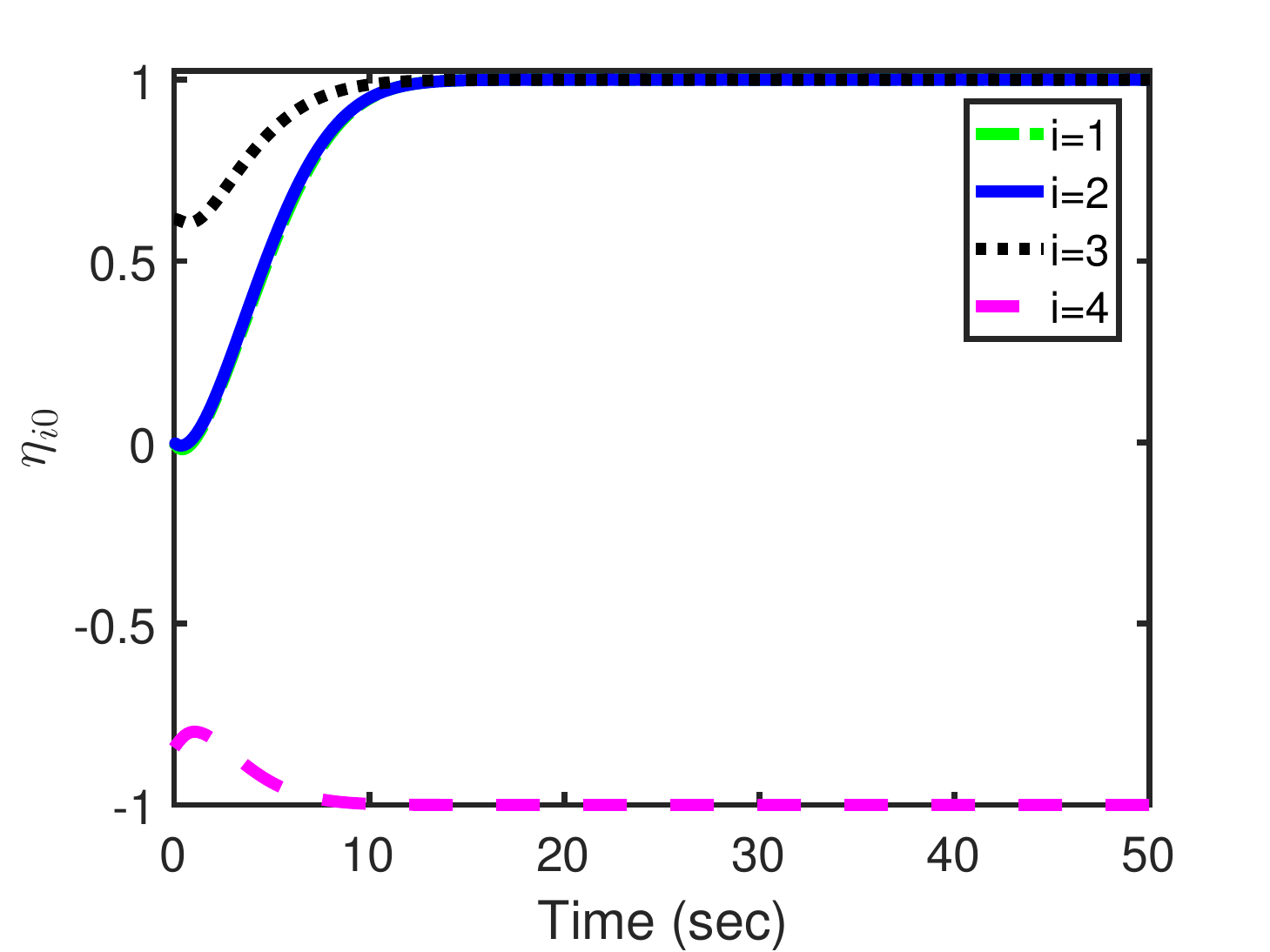} \label{fig3a}}
\subfloat[] {\includegraphics[height=4cm] {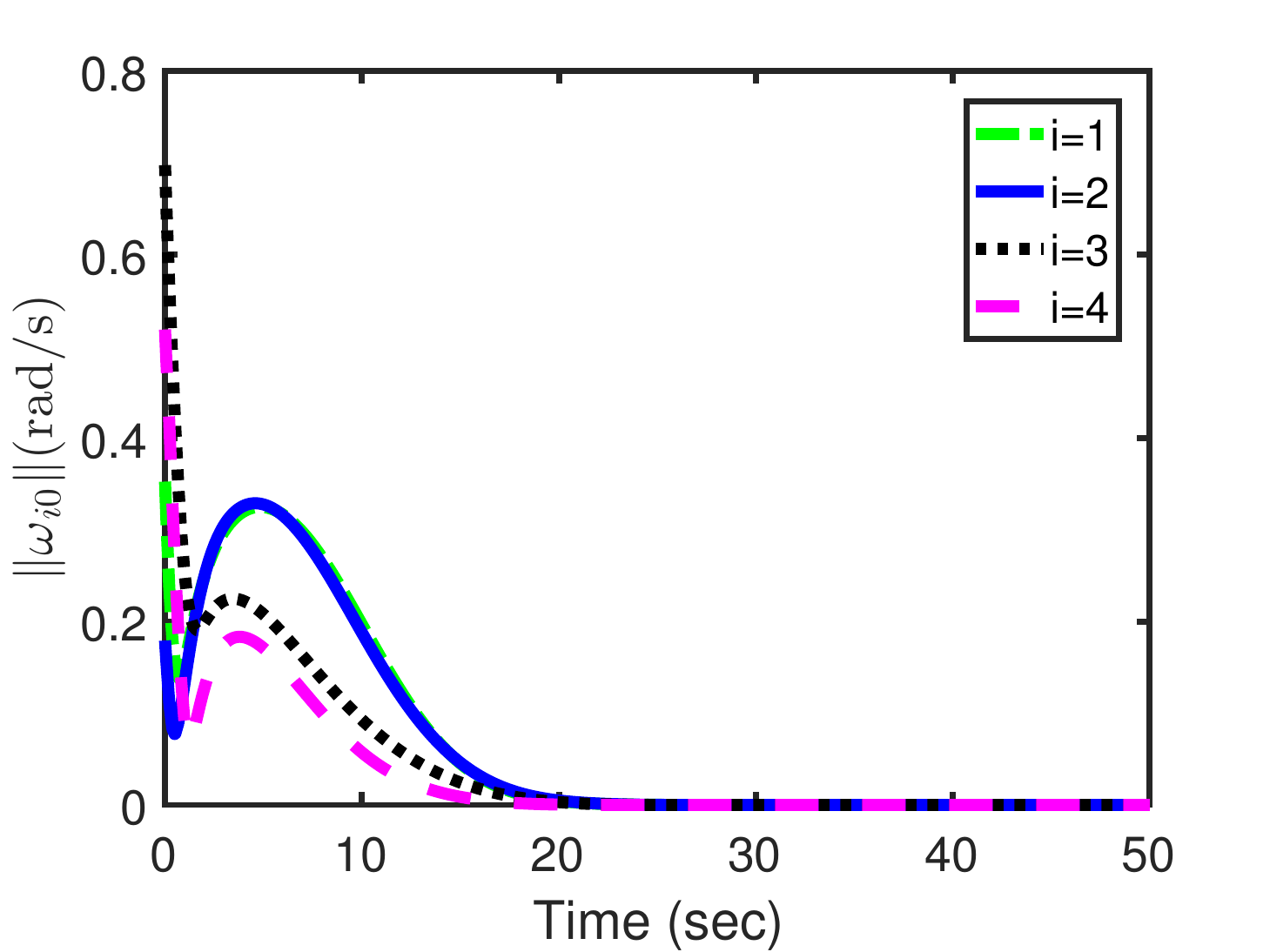} \label{fig3b}}
\subfloat[] {\includegraphics[height=4cm] {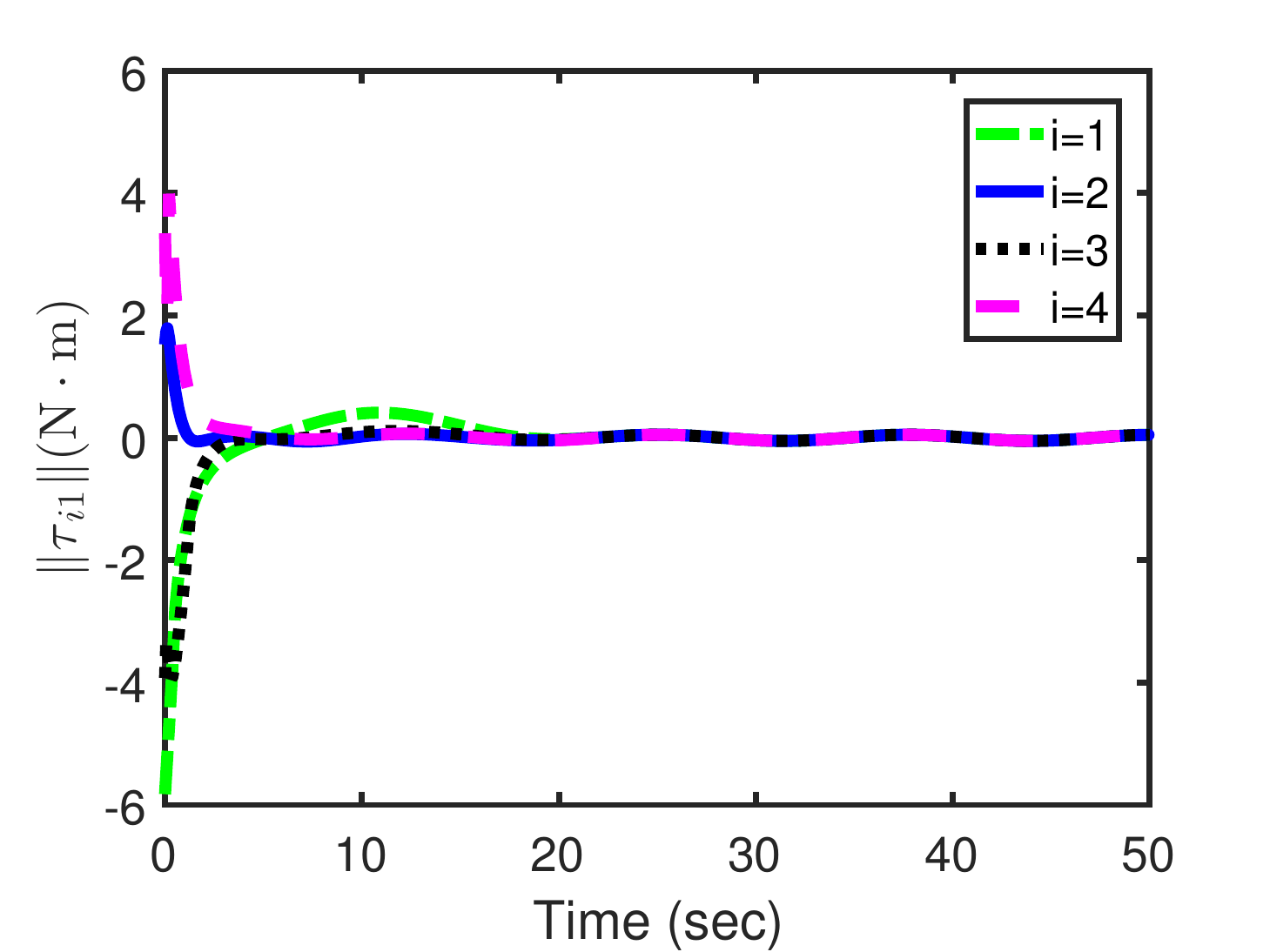} \label{fig3c}}
\caption{Simulation results using full-state feedback without any uncertainty: (a) $\eta_{i0}$, (b) $\|\tilde{\omega}_{i0}\|_2$, and (c) $\tau_{i1}$}  \label{fig3}
\end{center}
\end{figure*}

\begin{figure*}
\begin{center}
\subfloat[]{\includegraphics[height=4cm] {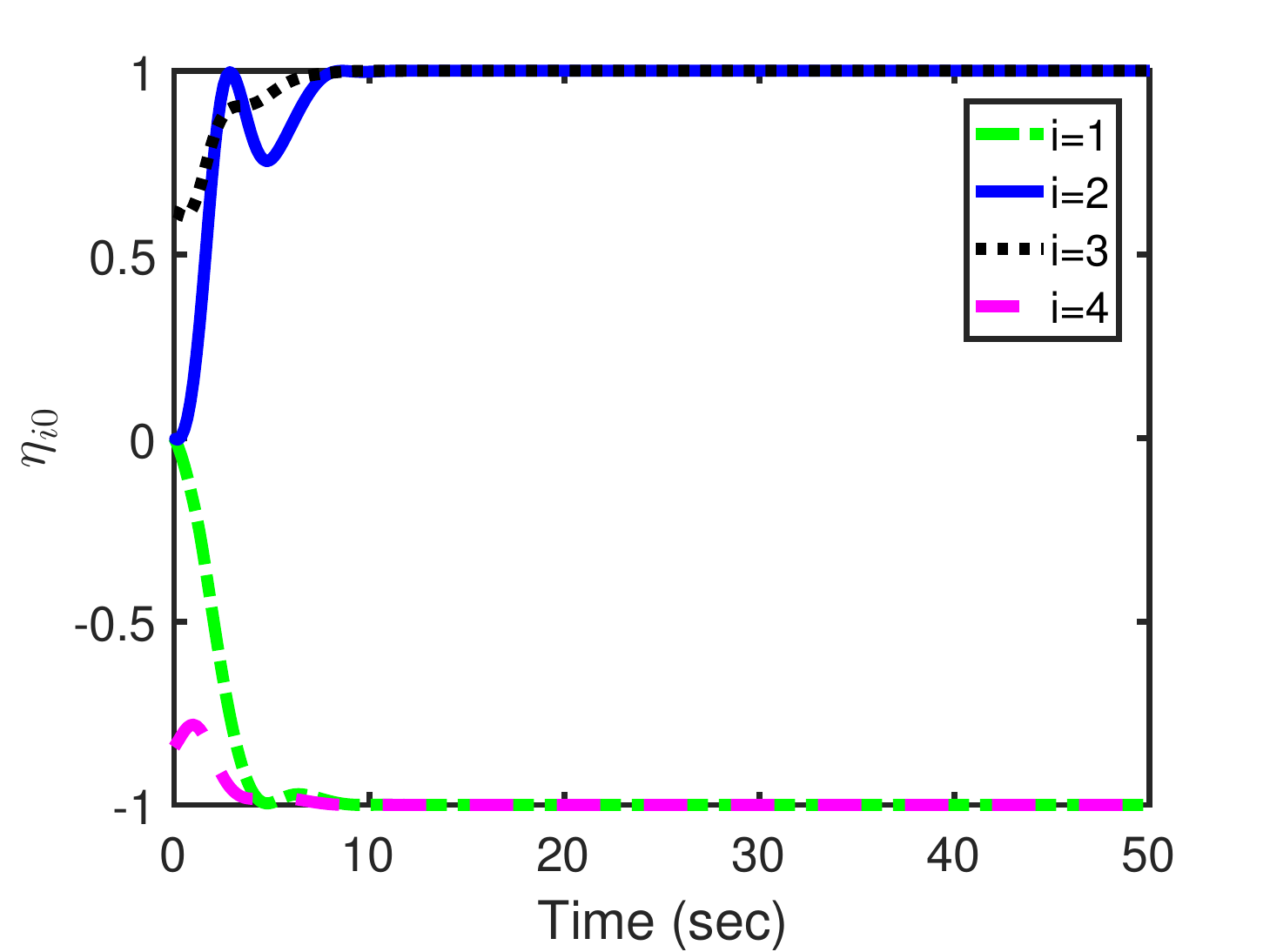} \label{fig4a}}
\subfloat[] {\includegraphics[height=4cm] {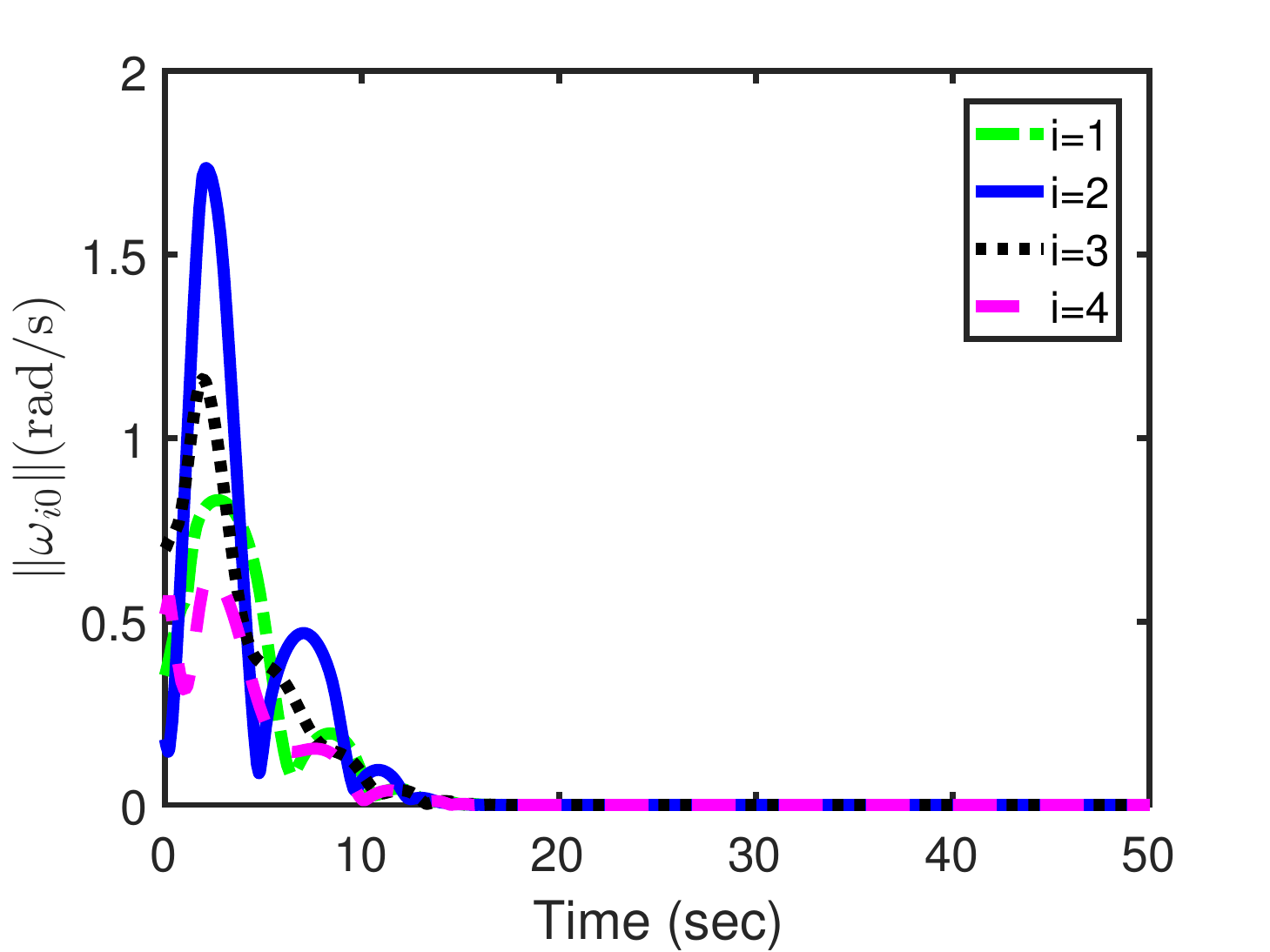} \label{fig4b}}
\subfloat[] {\includegraphics[height=4cm] {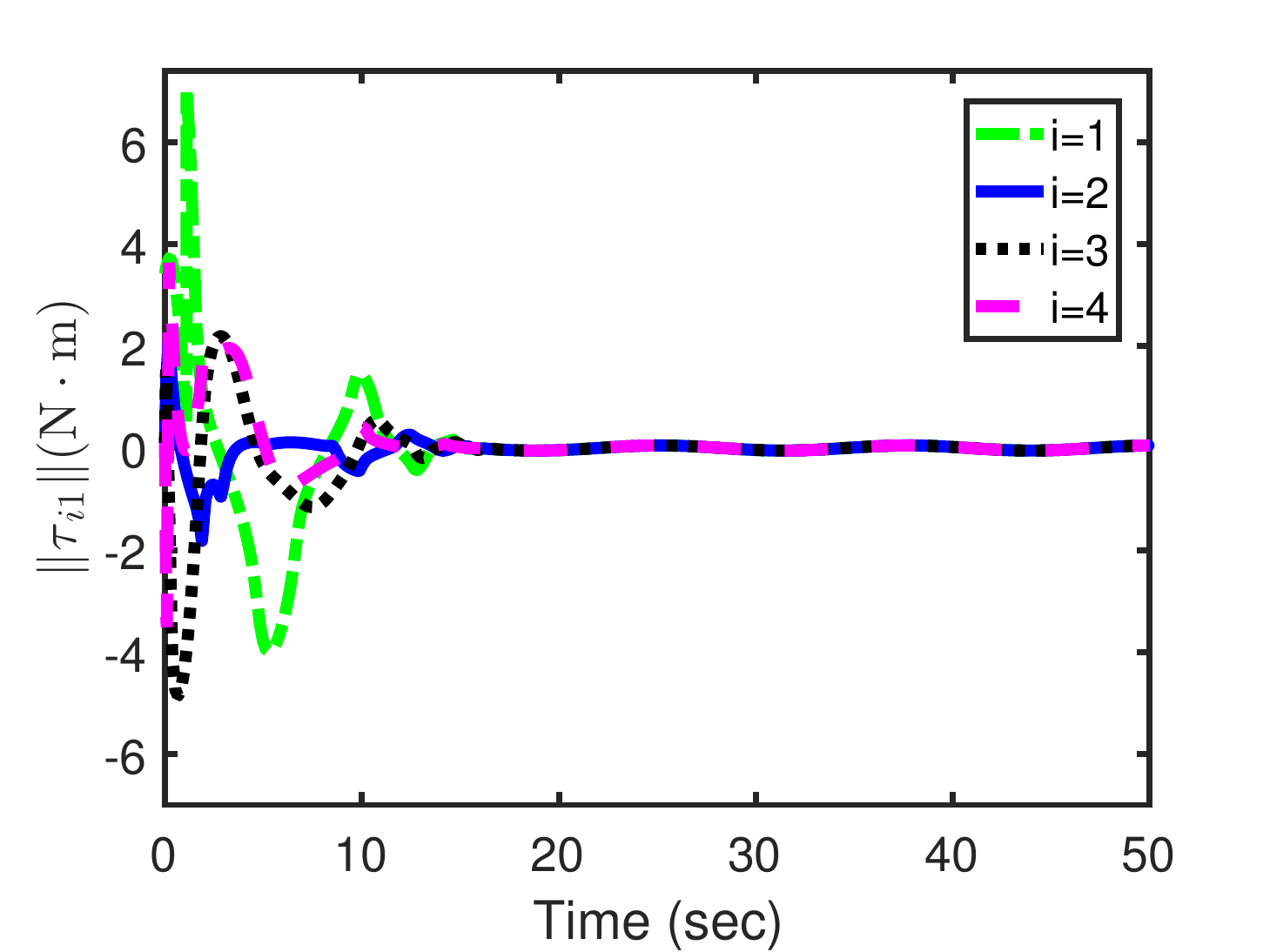} \label{fig4c}}
\caption{Simulation results using attitude-only feedback without any uncertainty: (a) $\eta_{i0}$, (b) $\|\tilde{\omega}_{i0}\|_2$, and (c) $\tau_{i1}$}  \label{fig4}
\end{center}
\end{figure*}

The first set of simulations assume that the attitude and angular velocity of each follower are available. Controller (\ref{eq_GFTC2}) and the distributed observer given in (\ref{eq_obs1})-(\ref{eq_obs4}) are applied. The corresponding control parameters are summarized in Table~\ref{tab1}. Figure~\ref{fig2} shows the estimation errors of the leader's quaternion, angular velocity, and acceleration in terms of their respective 2-norms. As shown by Fig.~\ref{fig2}, the proposed distributed observer exhibits fast transient and each follower recovers the leader's attitude, angular velocity, and acceleration within 5 s. Since both $P_0$ and $-P_0$ represent the leader's attitude, Fig.~\ref{fig:unwinding} compares the responses of $\|P_i-P_0\|$ and $\|P_i + P_0\|$. At the beginning $P_4$ is much closer to $-P_0$ but the observer still stabilizes $P_4$ through a longer path to $P_0$ instead of $-P_0$. Therefore, the unwinding phenomenon occurred in the distributed observer. Fortunately, this problem does not degrade the performance of the consensus controllers due to their global stabilizing feature, as shown by Figs.~\ref{fig3} and \ref{fig4} later. Figure~\ref{fig3} depicts the time histories of the attitude tracking error $\eta_{i0}(t)$, 2-norm of the angular velocity tracking error $\|\omega_{i0}(t)\|_2$, and the control torque component $\tau_{i1}$. It can be seen that the proposed full-state feedback consensus scheme successfully aligns the attitude of each follower with the dynamic leader through the shortest path in finite-time.

Next, the angular velocity measurements are removed from all followers and the velocity-free controller (\ref{eq_GFTC3}) is simulated under the same initial conditions. The control gains are shown in Table~\ref{tab1}. Figure~\ref{fig4} presents the responses of $\eta_{i0}(t)$, $\|\omega_{i0}(t)\|_2$ and $\tau_{i1}$. The followers also reach an agreement with the leader in finite time, though there is no angular velocity measurements. Compared to the full-sate feedback case, the angular velocity tracking error in Fig.~(\ref{fig4}) exhibits increased transients. This is mainly caused by the weakened damping effect due to the lack of velocity measurements for feedback control.



\subsection{Simulations with Time Delays and Disturbances}

In order to examine the robustness of the proposed methods, the measurement and communication updates of each spacecraft are delayed by 0.01 s. Additionally, we add a nonzero disturbance torque $d_i(t)=0.02[\cos(\theta_i t), \sin(\theta_i t), -\sin(\theta_i t)]^T \, \text{N}\cdot \text{m}$, where $\theta_i = 2\pi/(40+5i)$, $i\in\mathbb{I}_n$, to the $i$-th spacecraft. The control parameters remain the same as the previous example while the initial conditions of the four followers are renewed to

\begin{equation*}
\begin{matrix}
\begin{bmatrix}
  Q_1^T(0) \\
  Q_2^T(0) \\
  Q_3^T(0) \\
  Q_4^T(0)
\end{bmatrix} & = & \begin{bmatrix}
                      0.3162 & 0.1 & -0.8 & 0.5 \\
                      -0.1732 & -0.5 & 0.6 & -0.6 \\
                      0.7416 & 0.5 & -0.2 & -0.4 \\
                      -0.6557 & 0.5 & -0.4 & 0.4
                    \end{bmatrix}
\end{matrix}
\end{equation*}

\begin{equation*}
\begin{matrix}
\begin{bmatrix}
  \omega_1^T(0) \\
  \omega_2^T(0) \\
  \omega_3^T(0) \\
  \omega_4^T(0)
\end{bmatrix} & = & \begin{bmatrix}
                      0.2 & -0.1 & 0.4 \\
                      -0.5 & 0.6 & -0.6 \\
                      -0.5 & 0.4 & -0.2 \\
                      -0.1 & -0.6 & 0.1
                    \end{bmatrix}
\end{matrix} \, \textup{rad/s}
\end{equation*}

\begin{figure*}
\begin{center}
\subfloat[]{\includegraphics[height=4cm] {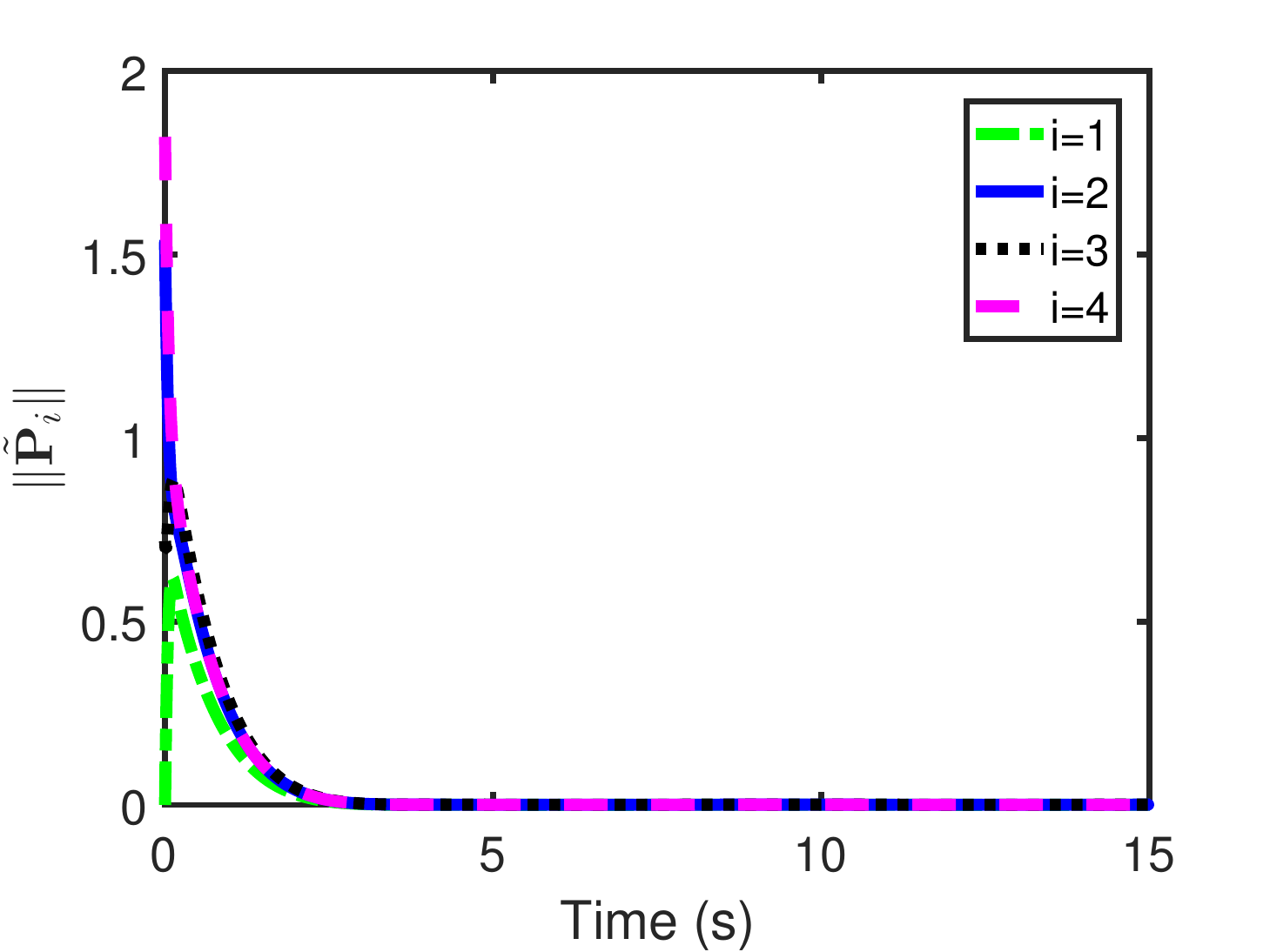} \label{fig5a}}
\subfloat[]{\includegraphics[height=4cm] {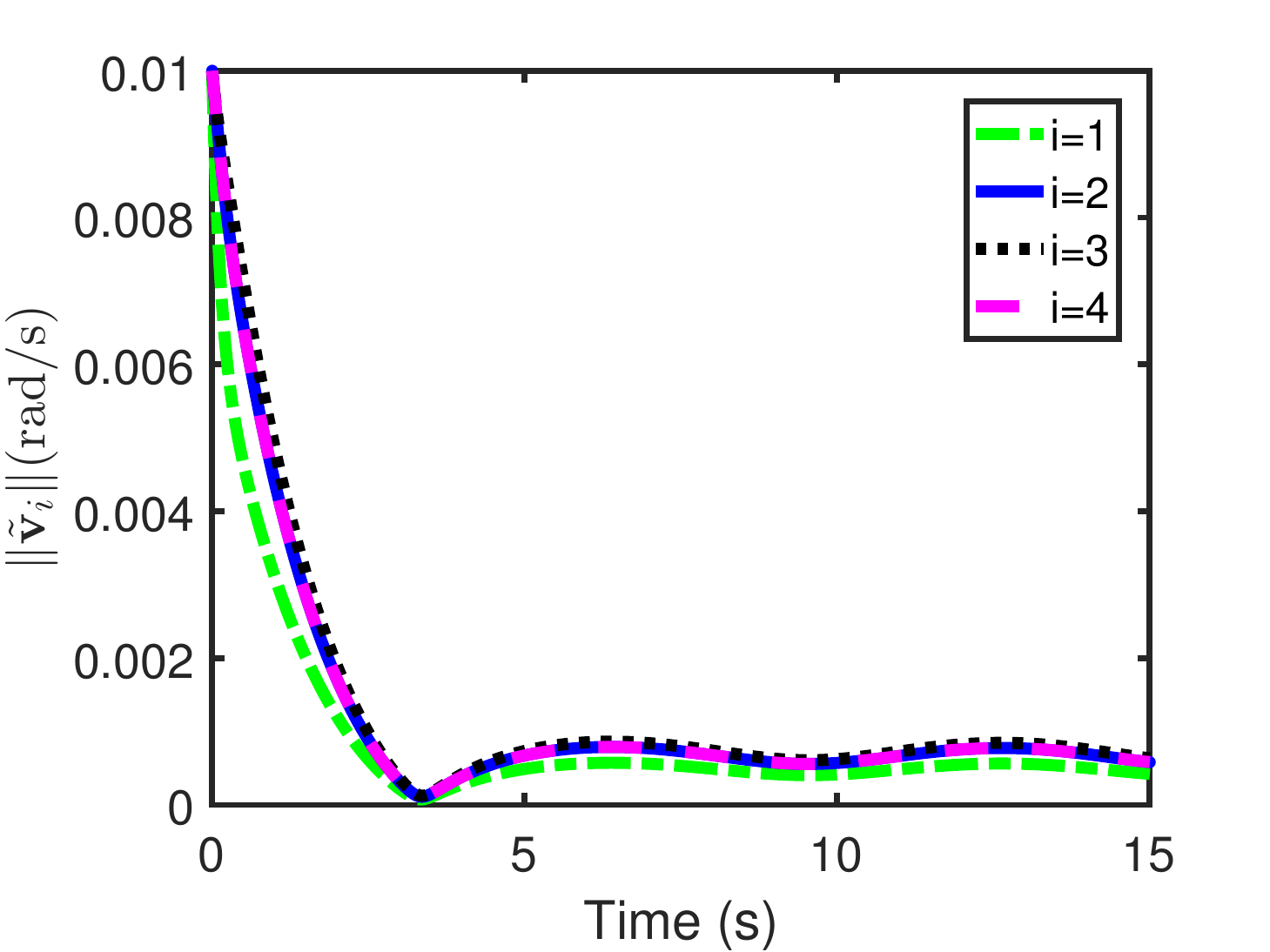} \label{fig5b}}
\subfloat[]{\includegraphics[height=4cm] {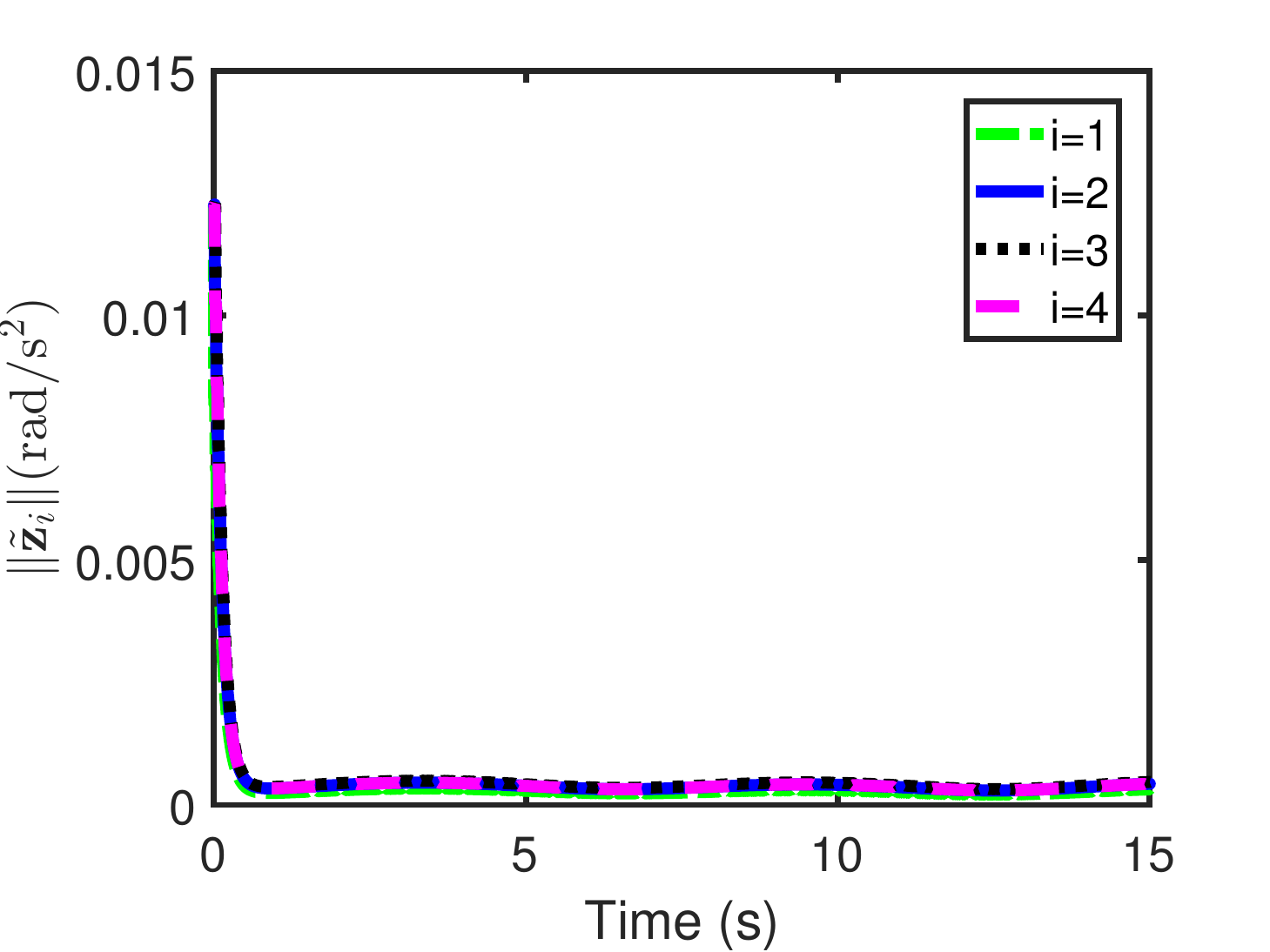} \label{fig5c}}
\caption{Response of the distributed observer with time delays and disturbances: (a) $\|\tilde{P}_i\|_2$, (b) $\|\tilde{v}_i\|_2$, and (c) $\|\tilde{z}_i\|_2$}  \label{fig5}
\end{center}
\end{figure*}

\begin{figure*}
\begin{center}
\subfloat[]{\includegraphics[height=4cm] {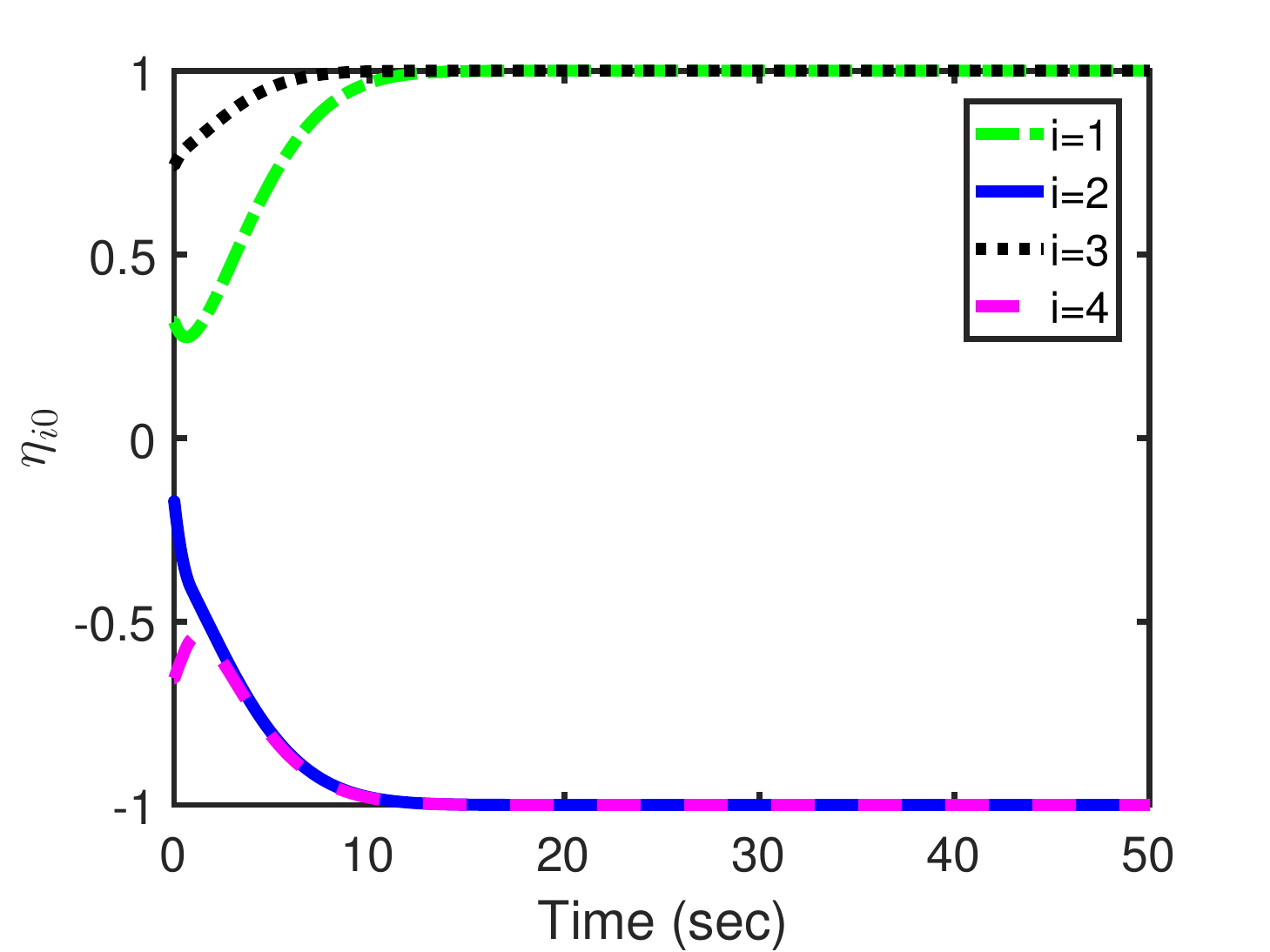} \label{fig6a}}
\subfloat[] {\includegraphics[height=4cm] {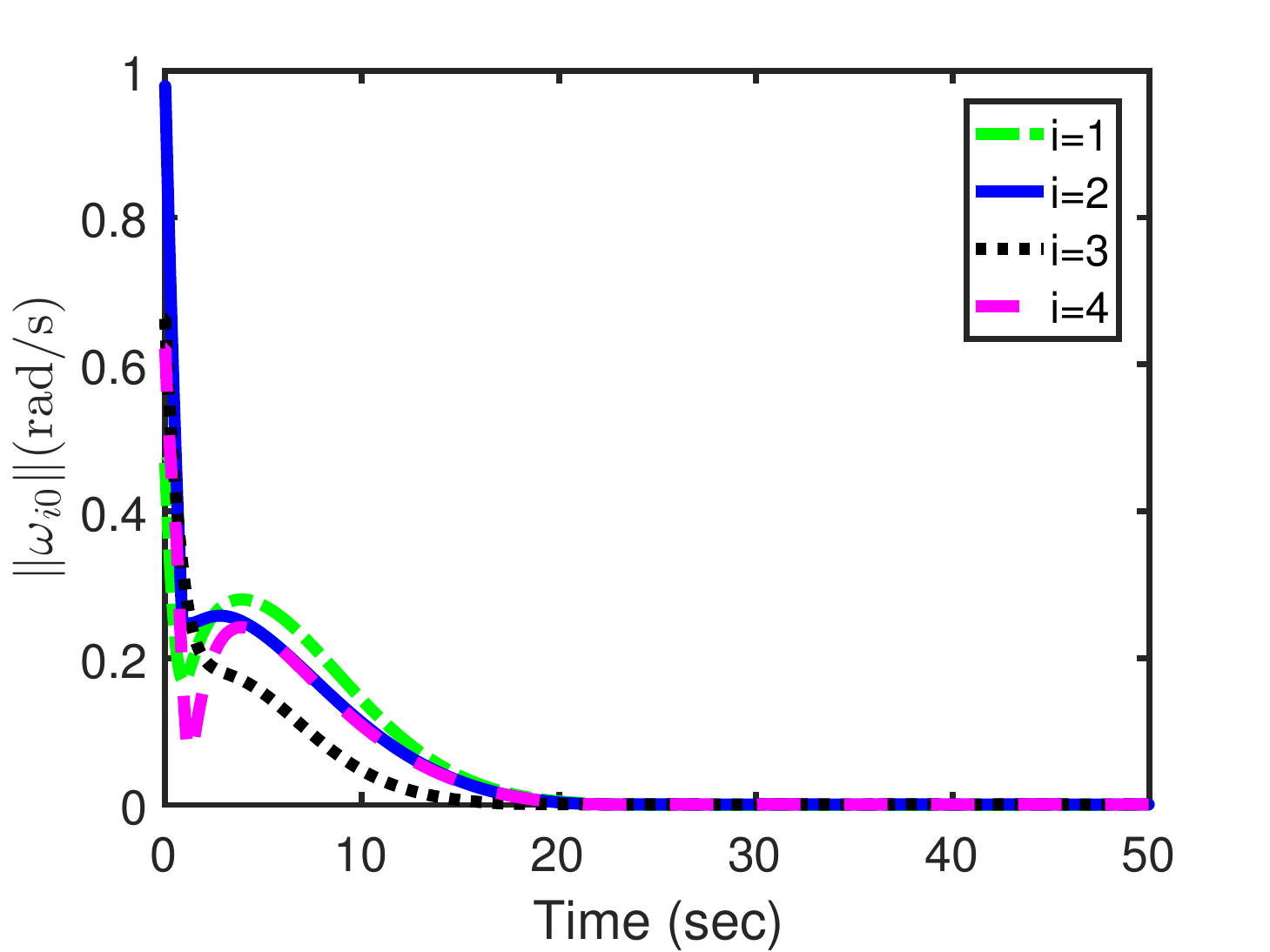} \label{fig6b}}
\subfloat[] {\includegraphics[height=4cm] {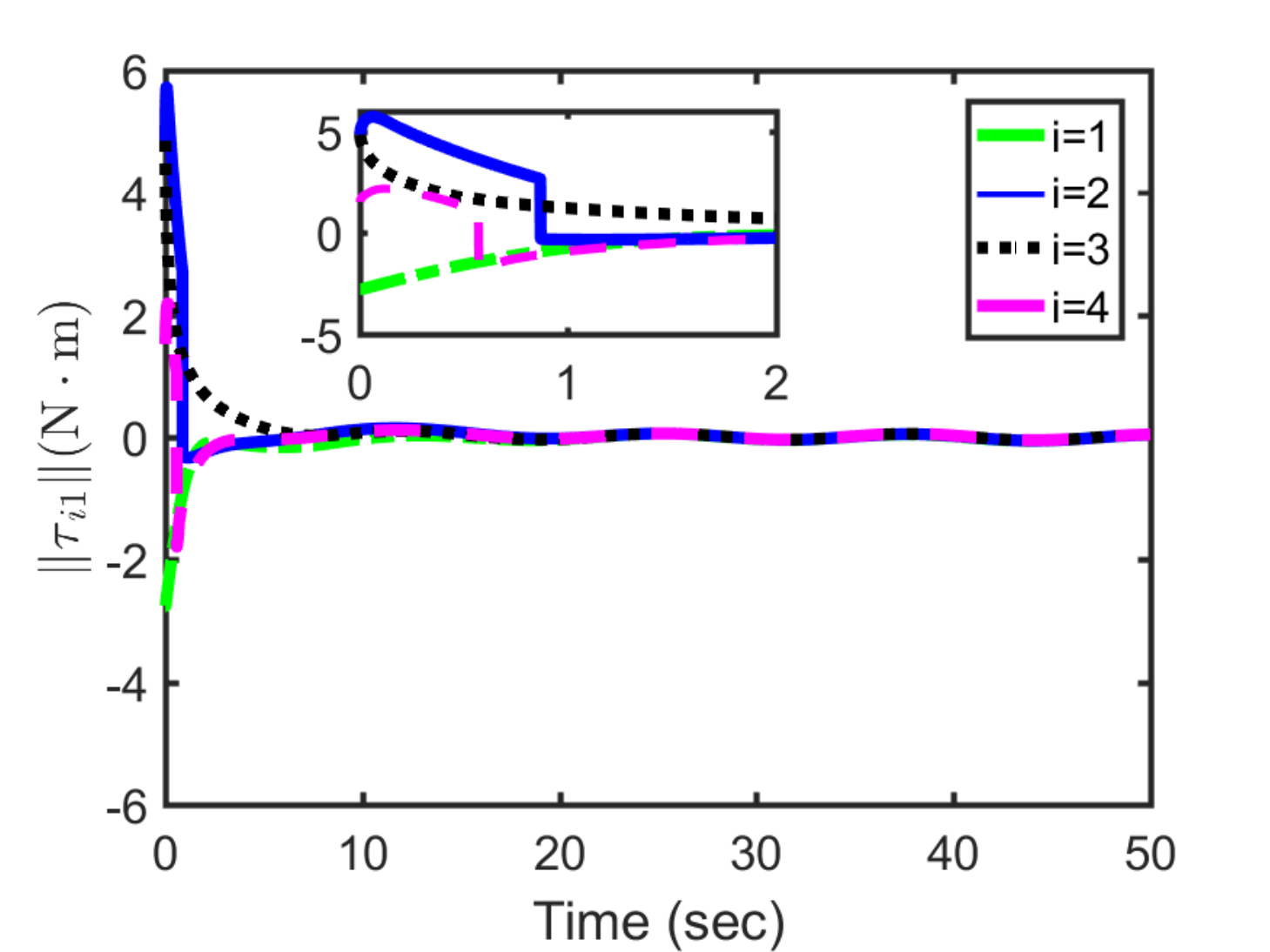} \label{fig6c}}
\caption{Simulation results using full-state feedback with time delays and disturbances: (a) $\eta_{i0}$, (b) $\|\tilde{\omega}_{i0}\|_2$, and (c) $\tau_{i1}$}  \label{fig6}
\end{center}
\end{figure*}

\begin{figure*}
\begin{center}
\subfloat[]{\includegraphics[height=4cm] {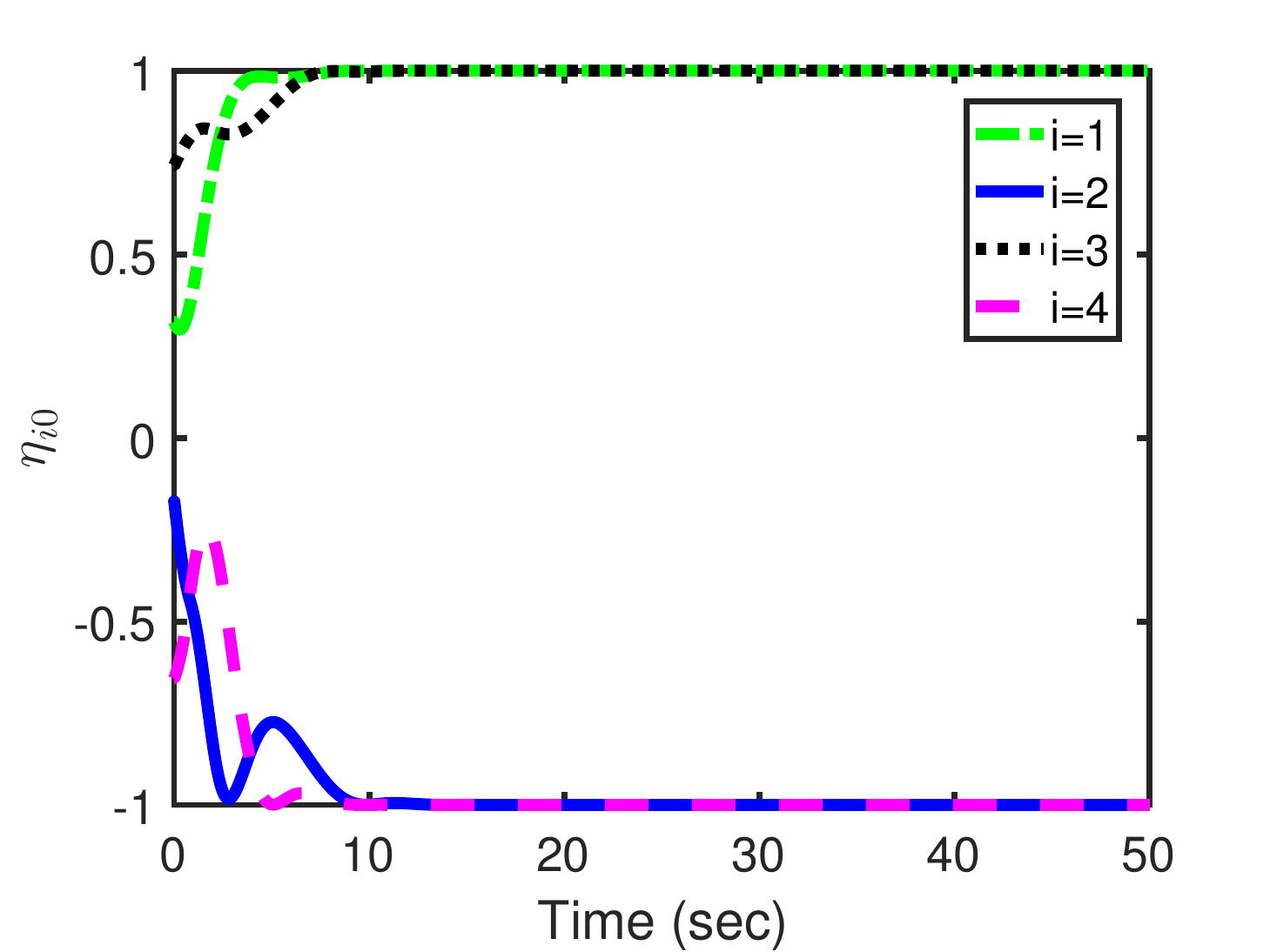} \label{fig7a}}
\subfloat[] {\includegraphics[height=4cm] {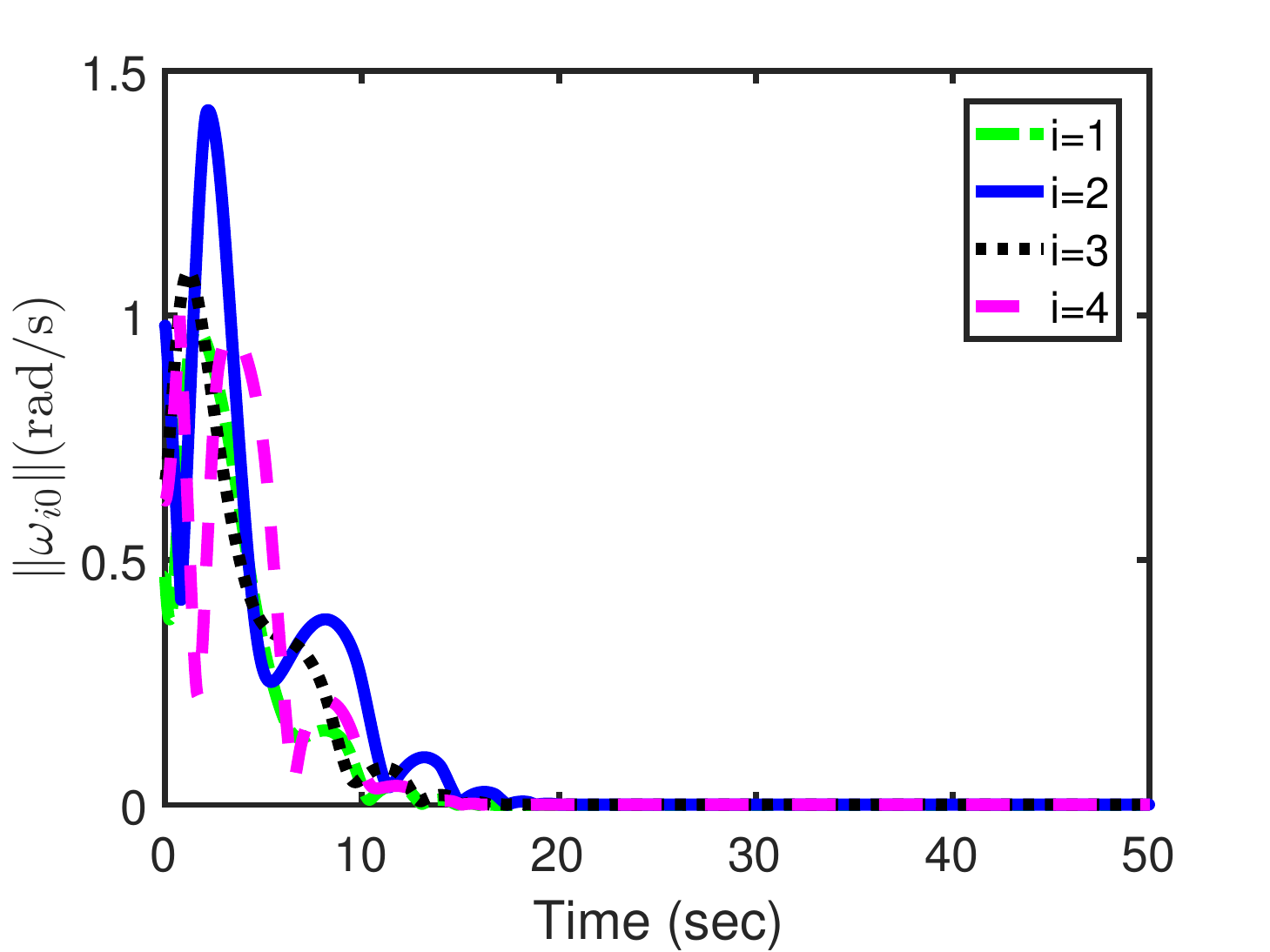} \label{fig7b}}
\subfloat[] {\includegraphics[height=4cm] {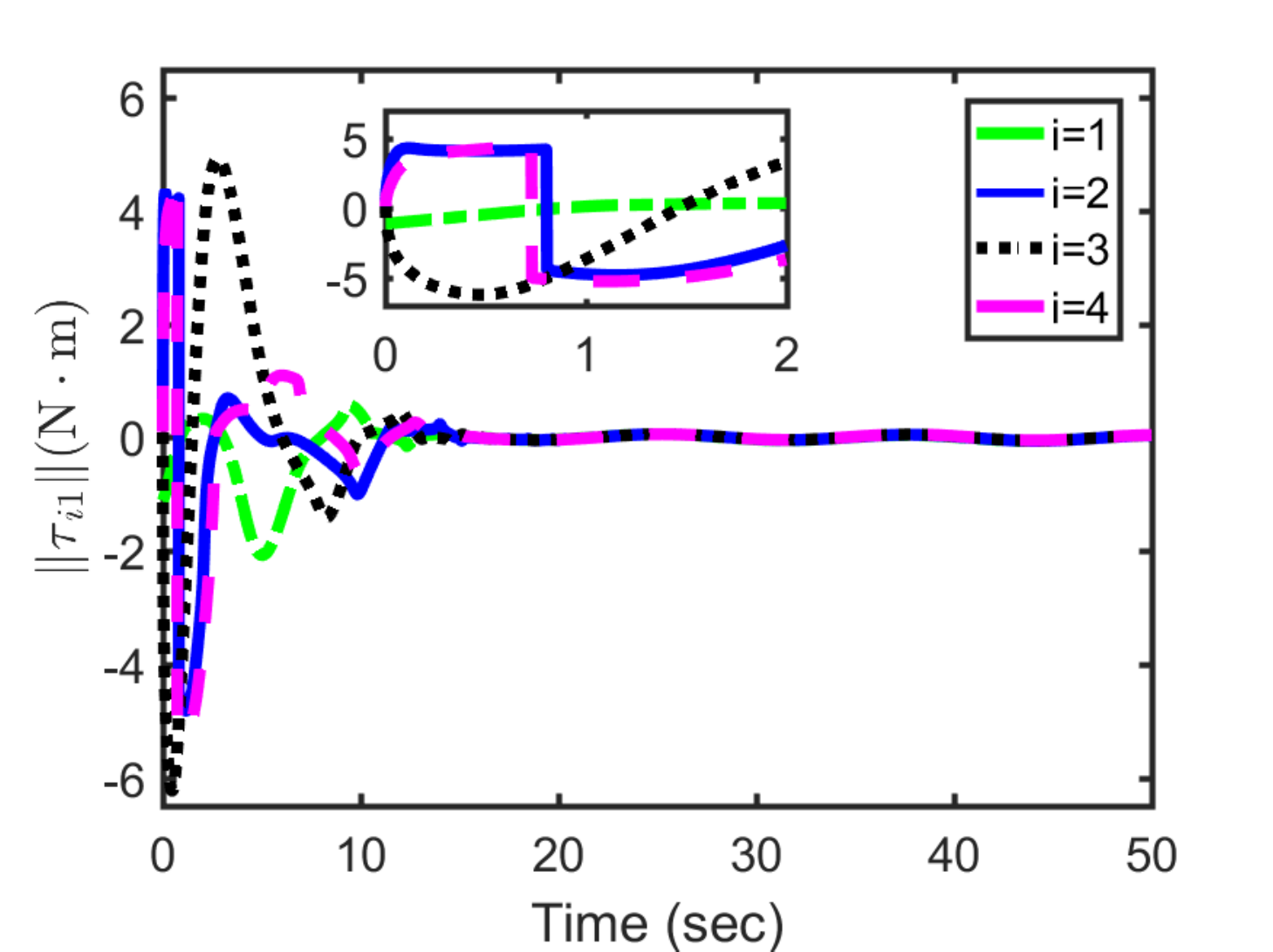} \label{fig7c}}
\caption{Simulation results using attitude-only feedback with time delays and disturbances: (a) $\eta_{i0}$, (b) $\|\tilde{\omega}_{i0}\|_2$, and (c) $\tau_{i1}$}  \label{fig7}
\end{center}
\end{figure*}

The simulation results of the distributed observer and full-state feedback controller are presented in Figures~\ref{fig5} and \ref{fig6}. In particular, Fig.~\ref{fig5} shows the estimation errors of the leader's quaternion, angular velocity, and acceleration in terms of their respective 2-norms while Fig.~\ref{fig6} plots the time histories of $\eta_{i0}(t)$, $\|\omega_{i0}(t)\|_2$, and $\tau_{i1}$. In spite of the time delays and external disturbances, the proposed distributed observer still maintains fast transient and each follower obtains the leader's trajectory with high accuracy within 5 s. The full-state feedback controller successfully synchronizes the attitude of each follower with the dynamic leader (Figs.~\ref{fig6a} and \ref{fig6b}). From Fig.~\ref{fig6c}, it can be observed that jumps occur in the control torques of followers 2 and 4 as time approaches 1 s while the control torques of followers 1 and 3 remain continuous during the entire control phase. Figure~\ref{fig7} presents the responses of $\eta_{i0}(t)$, $\|\omega_{i0}(t)\|_2$ and $\tau_{i1}$ for the attitude-only feedback controller. Similarly to the full-state feedback case, the control torques of followers 2 and 4 switched once before 1 s to change the rotation direction (Fig.~\ref{fig7c}). By means of controller~(\ref{eq_GFTC3}), the followers reach an agreement with the leader while avoiding the unwinding phenomenon, despite the absence of angular velocity feedback (Figs.~\ref{fig7a} and \ref{fig7b}). The simulation results shows that the proposed methods possess some robustness to small time delays and disturbances.


\subsection{Simulations with Time-Varying Topologies}

\begin{figure}
\begin{center}
\includegraphics[height=4cm]{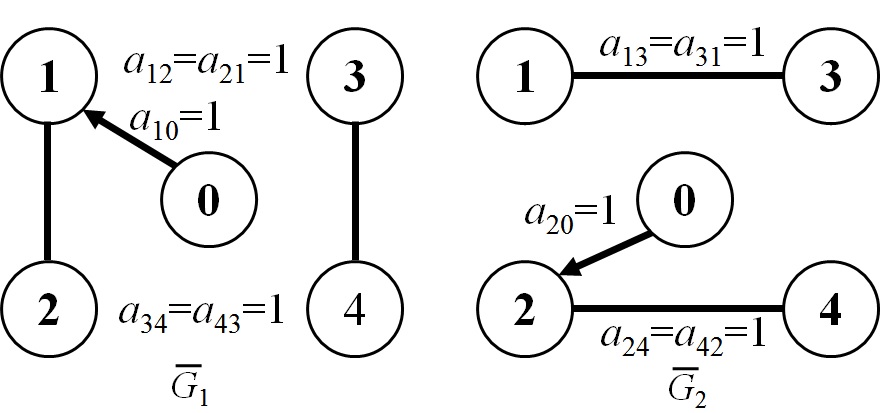}
\caption{Communication graph of the leader-following system}
\label{fig8}
\end{center}
\end{figure}

Although the design in this paper and the previous two examples both assume fixed communication topologies among the multi-spacecraft system, it is still interesting the see how the proposed methods behave with time-varying communication topologies. An example is shown in Fig.~\ref{fig8}, the graph $\bar{G}(t)$ has two possible topologies $\bar{G}_1$ and $\bar{G}_2$. Let $\bar{G}(t)=\bar{G}_1$ for $t\in[m t_D, m t_D + 0.5t_D)$ and $\bar{G}(t)=\bar{G}_2$ for $t\in[m t_D + 0.5t_D, (m+1) t_D)$, where $t_D = 0.2$ s and $m=0,1,2,\cdots$. Clearly, the graph $\bar{G}(t)$ is not connected at any moment but the union graph $\bigcup_{t\in[t,t+t_D)} \bar{G}(t)$ is connected for any $t\geq 0$. Hence, $\bar{G}(t)$ satisfies the joint strong connectivity \cite{Thunberg:14}, a condition that is usually employed to guarantee consensus in multi-agent systems with time-varying topologies.

\begin{figure*}
\begin{center}
\subfloat[]{\includegraphics[height=4cm] {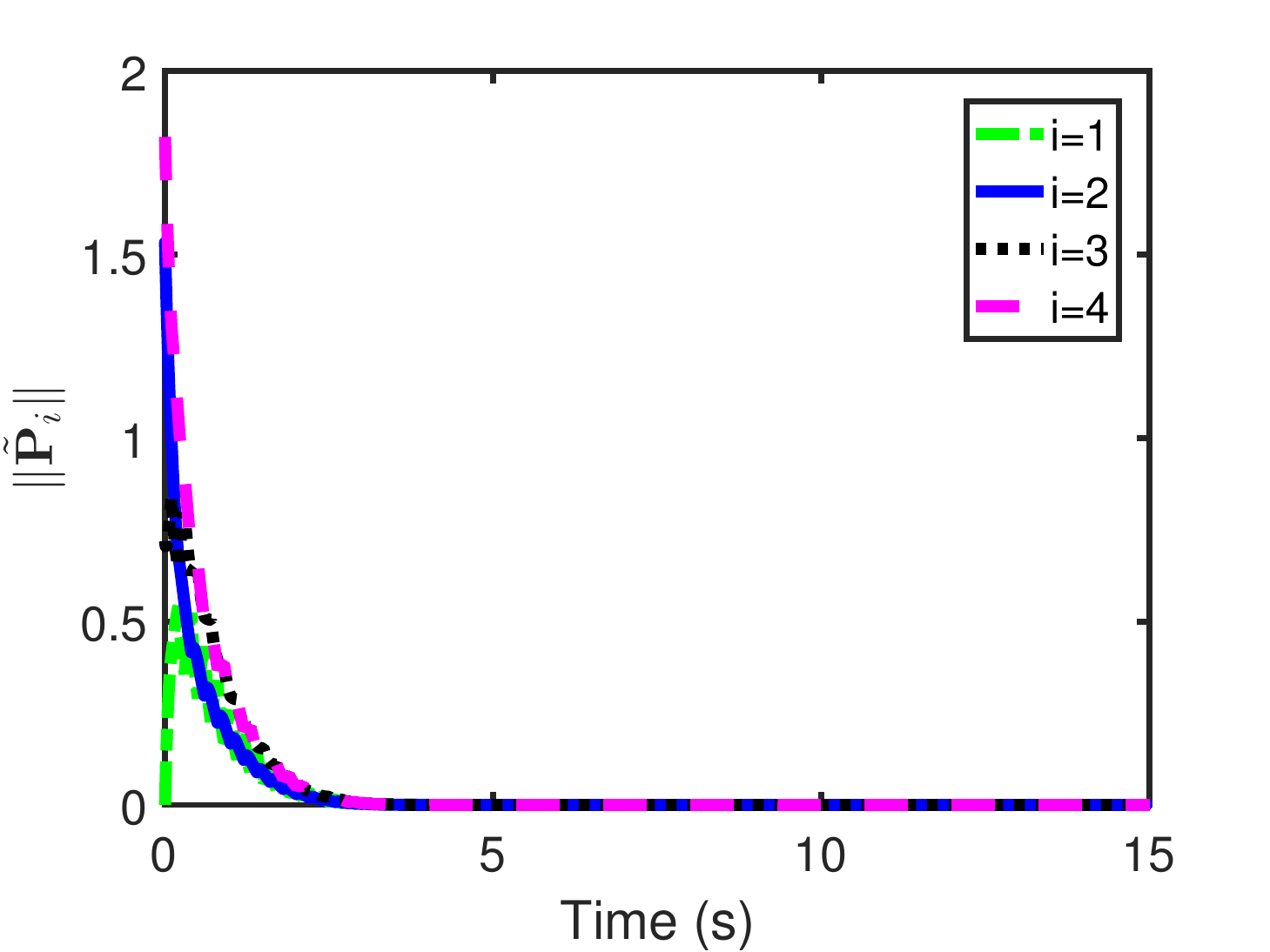} \label{fig9a}}
\subfloat[]{\includegraphics[height=4cm] {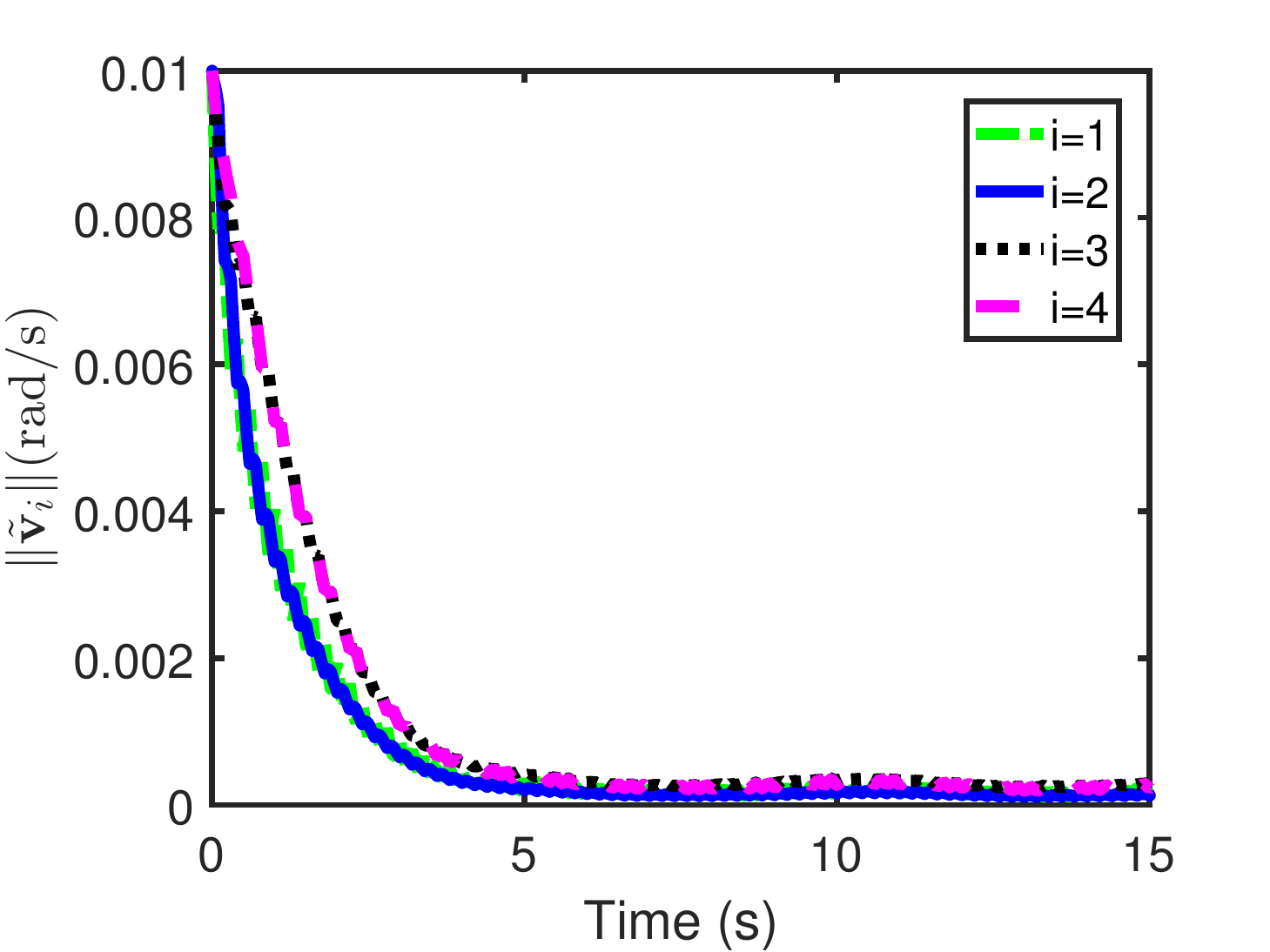} \label{fig9b}}
\subfloat[]{\includegraphics[height=4cm] {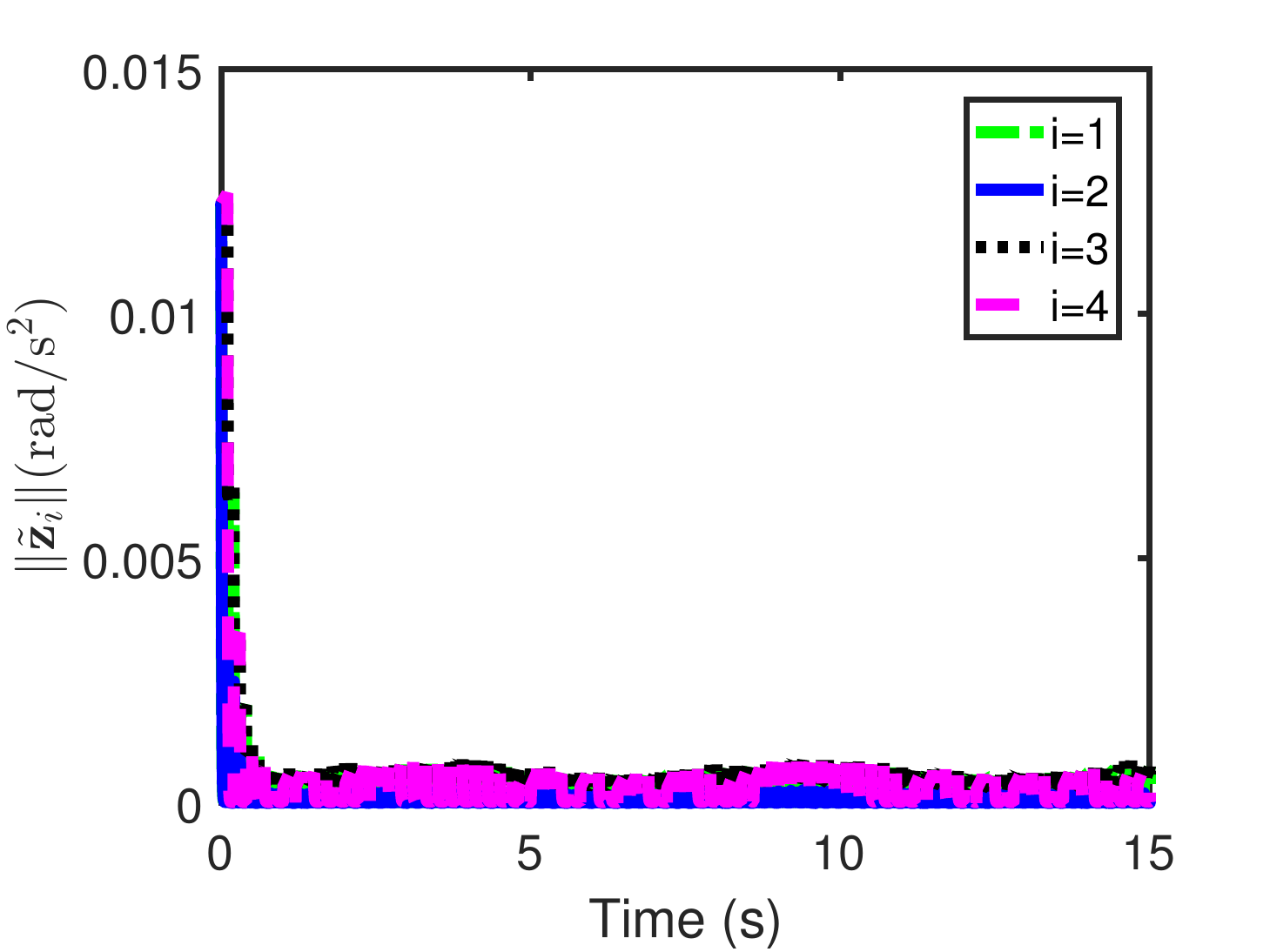} \label{fig9c}}
\caption{Response of the distributed observer with time-varying topologies: (a) $\|\tilde{P}_i\|_2$, (b) $\|\tilde{v}_i\|_2$, and (c) $\|\tilde{z}_i\|_2$}  \label{fig9}
\end{center}
\end{figure*}

\begin{figure*}
\begin{center}
\subfloat[]{\includegraphics[height=4cm] {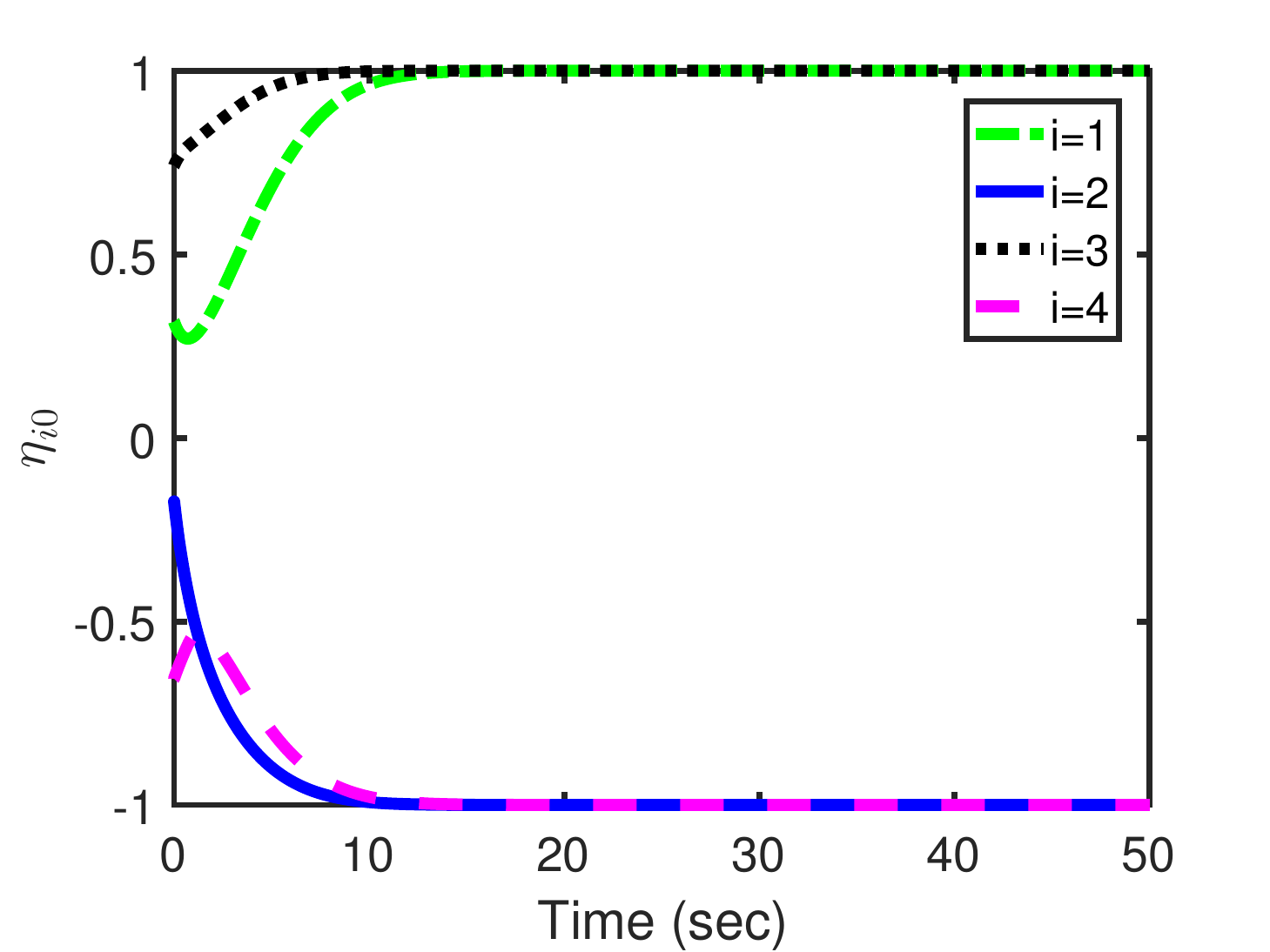} \label{fig10a}}
\subfloat[] {\includegraphics[height=4cm] {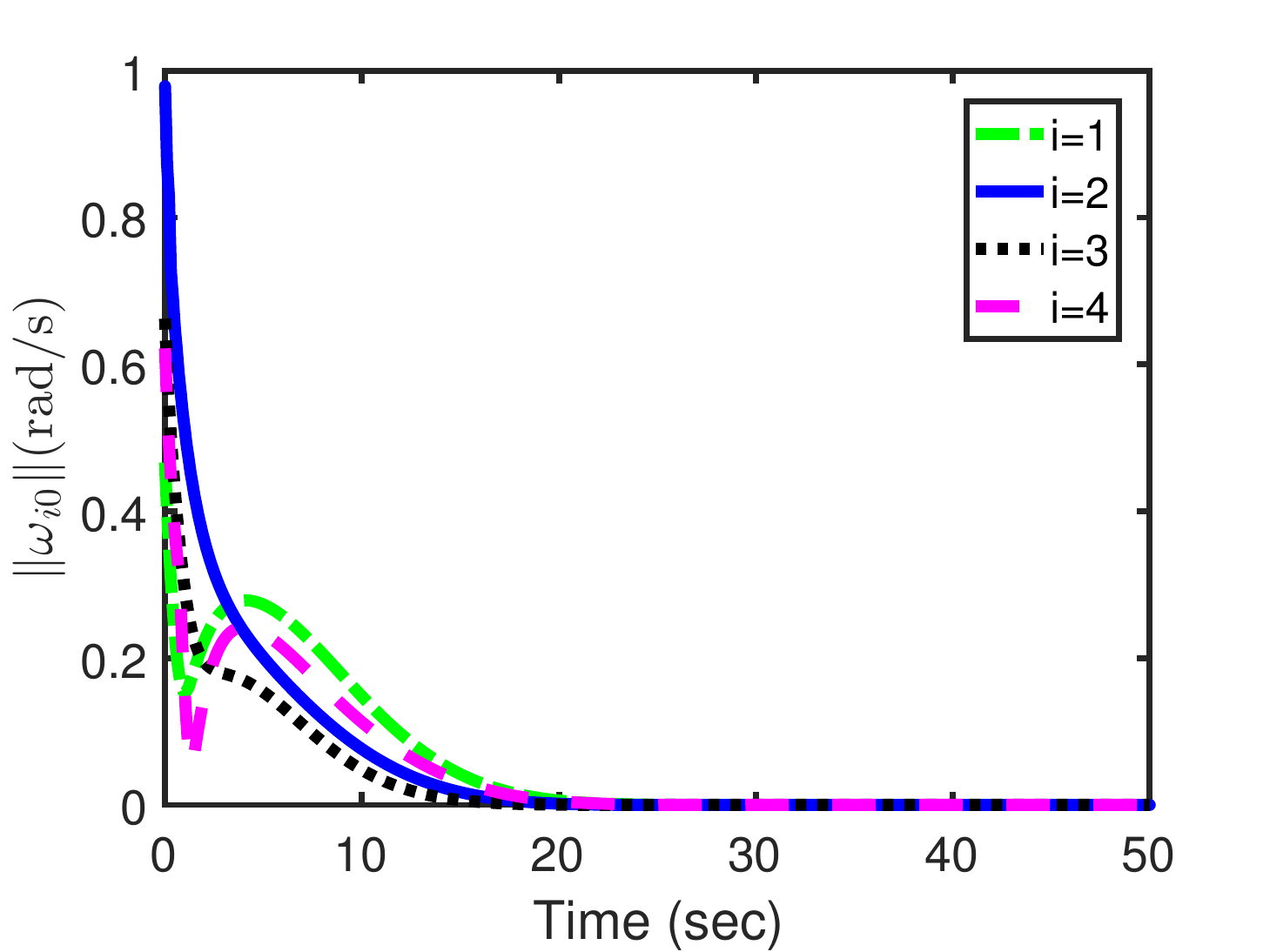} \label{fig10b}}
\subfloat[] {\includegraphics[height=4cm] {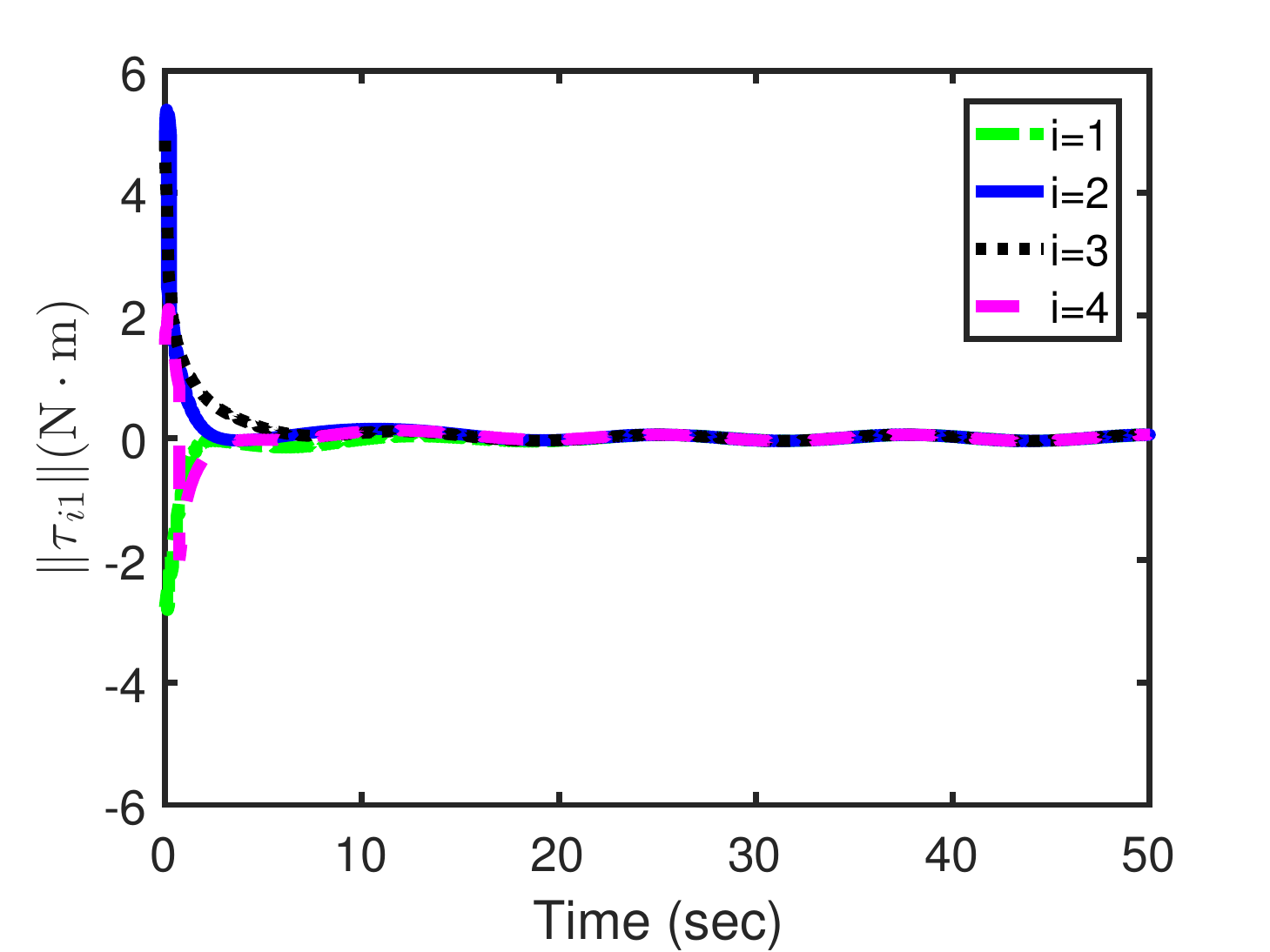} \label{fig10c}}
\caption{Simulation results using full-state feedback with time-varying topologies: (a) $\eta_{i0}$, (b) $\|\tilde{\omega}_{i0}\|_2$, and (c) $\tau_{i1}$}  \label{fig10}
\end{center}
\end{figure*}

\begin{figure*}
\begin{center}
\subfloat[]{\includegraphics[height=4cm] {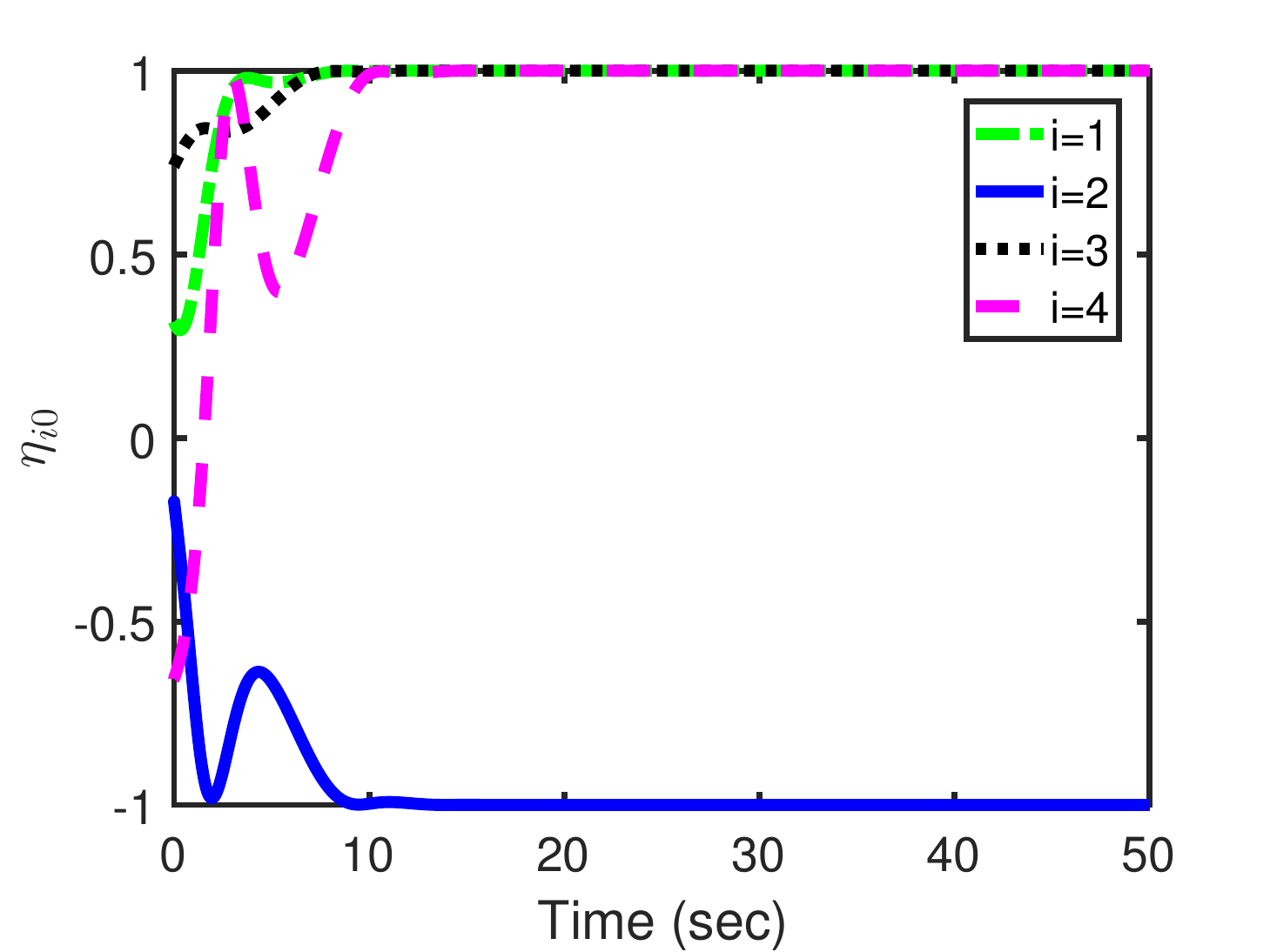} \label{fig11a}}
\subfloat[] {\includegraphics[height=4cm] {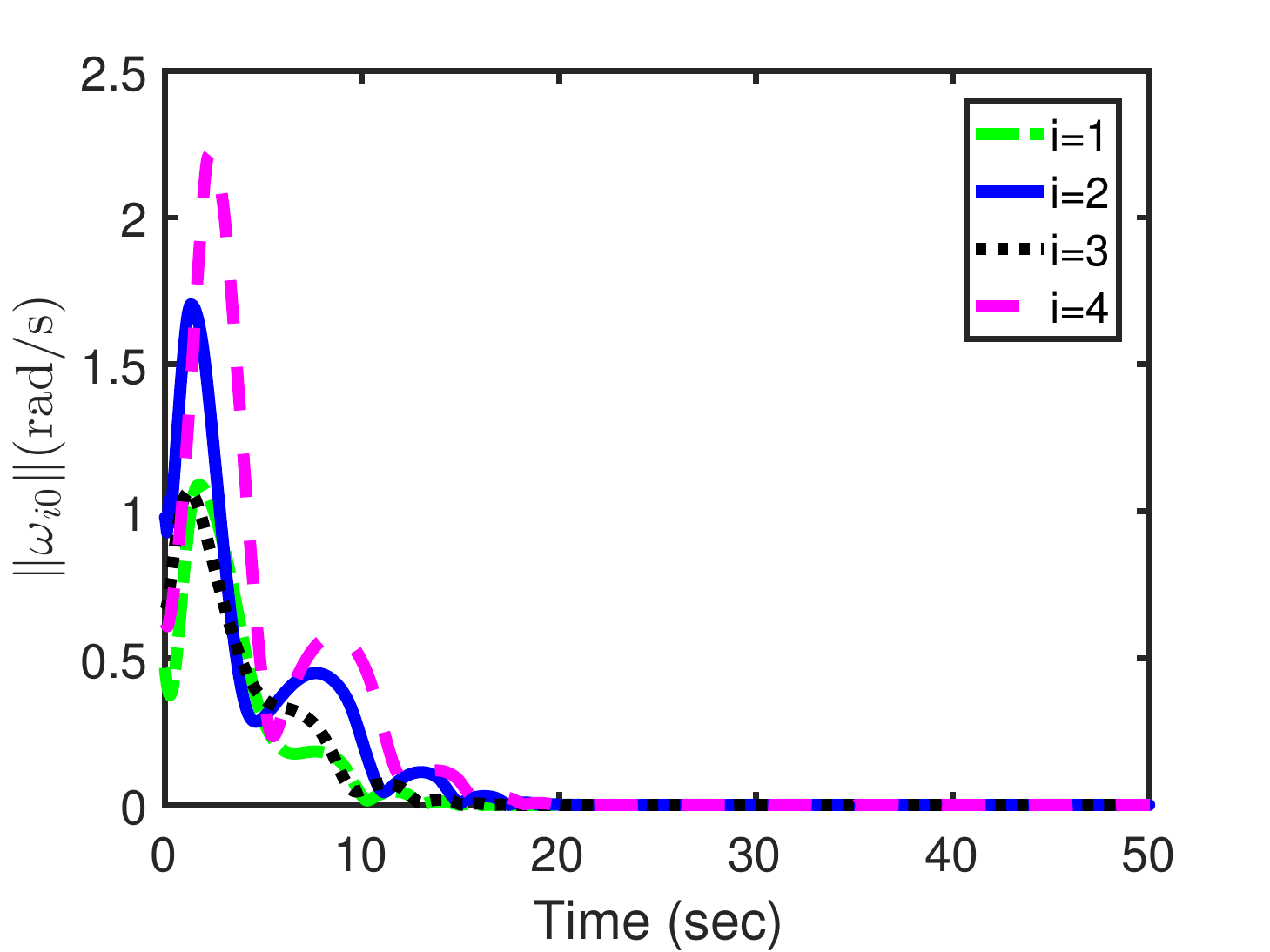} \label{fig11b}}
\subfloat[] {\includegraphics[height=4cm] {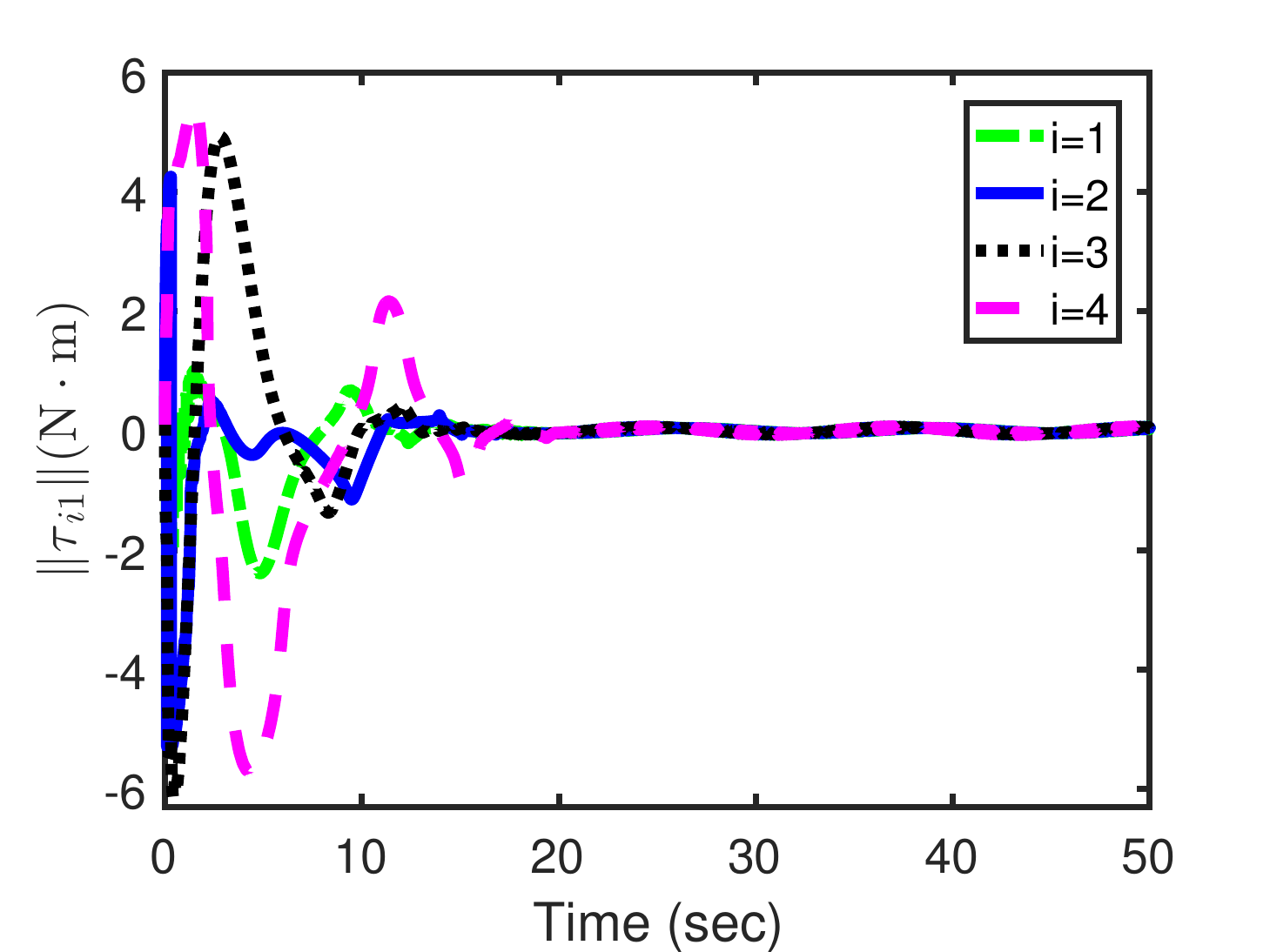} \label{fig11c}}
\caption{Simulation results using attitude-only feedback with time-varying topologies: (a) $\eta_{i0}$, (b) $\|\tilde{\omega}_{i0}\|_2$, and (c) $\tau_{i1}$}  \label{fig11}
\end{center}
\end{figure*}

In order to clearly see the effect of the switching topology, simulations are conducted without external disturbance and time delays in measurements and communications. The observer and controller gains are given in Table~\ref{tab1} while the initial attitudes and angular velocities remain the same as those in Section V.B. The simulation results are presented in Figs.~\ref{fig9}--\ref{fig11}, showing the responses of the distributed observer, the full-state feedback controller, and the attitude-only controller, respectively. It can be seen that the estimation errors and attitude tracking errors are still convergent with the considered switching topologies. The numerical results inspire us to conjecture that the proposed distributed observer and controller also applies to time-varying topologies satisfying the joint strong connectivity. A rigorous proof of this conjecture is left for future research.


\section{Conclusions}
Quaternion-based attitude consensus schemes were proposed for a group of leader-following spacecraft with an undirected and connected communication graph among the followers. Instrumental in our approach is a nonlinear distributed observer for the leader's states, which establishes second-order sliding modes for attitude and angular velocity estimation errors and hence recovers the leader's trajectory in finite time for each follower. By appropriately integrating the observer and two quaternion-based hybrid homogeneous controllers originally developed for single spacecraft, global group attitude agreement was obtained in finite time respectively with full-state measurements and attitude-only measurements. It is worth noting that the proposed distributed observer allows any reference trajectories with bounded time derivatives and can also be combined with many other single-spacecraft attitude controllers to achieve cooperative attitude tracking while pursuing different performance.
\bibliographystyle{IEEEtran}
\bibliography{ref_gui}

\appendices
\section{Proof of Lemmas \ref{lem_E1} and \ref{lem_E2}}
\label{app1}       
\textbf{Proof of Lemma \ref{lem_E1}}. The first identity of Lemma \ref{lem_E1} can be verified by direct computations and hence is omitted. According the definition of $E(\cdot)$ in (\ref{eq_kine}), it follows that
$$ E^{T}(Q^{\prime}) Q = -\eta q^{\prime} + \eta^{\prime}q - q^{\prime}\times q = -\textup{vec}(Q^*\circ{Q^\prime}),$$
which is the transpose of the second identity of Lemma \ref{lem_E1}. As a special case, letting $Q=Q^{\prime}$ yields ${Q^T}E(Q)=0$.

\textbf{Proof of Lemma \ref{lem_E2}}. Next, applying Lemma \ref{lem_E1} one can derive
\begin{equation} \label{eq_A1}
\begin{split}
\sum_{i=1}^{n}\sum_{j=1}^{n} b_{ij}{Q_{j}^T}E(Q_i)&= -\sum_{i=1}^{n}\sum_{j=1}^{n} b_{ij} \textup{vec}^{T}(Q_{j}^*\circ{Q_i}) \\
&= -\sum_{j=1}^{n}\sum_{i=1}^{n} b_{ji} \textup{vec}^{T}(Q_{i}^*\circ{Q_j}).
\end{split}
\end{equation}
Invoking $b_{ij}=b_{ji}$ and $-\textup{vec}^{T}(Q_{i}^*\circ{Q_j})= \textup{vec}^{T}(Q_{j}^*\circ{Q_i})$, it follows from (\ref{eq_A1}) that
\begin{align*}
\sum_{i=1}^{n}\sum_{j=1}^{n} b_{ij}{Q_{j}^T}E(Q_i)&= \sum_{i=1}^{n}\sum_{j=1}^{n} b_{ij} \textup{vec}^{T}(Q_{j}^*\circ{Q_i}) \\
&= -\sum_{i=1}^{n}\sum_{j=1}^{n} b_{ij}{Q_{j}^T}E(Q_i),
\end{align*}
which implies $\sum_{i=1}^{n}\sum_{j=1}^{n} b_{ij}{Q_{j}^T}E(Q_i)=0$.

\section{Proof of Theorem \ref{thm_obs}}         
\label{app2}
The proof is divided into four steps by showing the convergence of the sliding mode differentiator and $z_i$, $v_i$, and $P_i$ successively. \\
\textit{Step 1: the convergence of (\ref{eq_errobs4})}. When $a_{i0}=0$ (i.e., no direct access to the leader), it follows from~(\ref{eq_obs4}) that $y_i(t)=w_i(t)=0$ and $\tilde{w}_i = -\dot{\omega}_0$. When $a_{i0}=1$ (i.e., direct access to the leader), applying Lemma~\ref{lem_STA} to~(\ref{eq_errobs4}) implies that $w_i(t)\to z_0(t)=\dot{\omega}_0(t)$ in a finite time $t_r$, i.e, $A_0\tilde{w}(t) = 0$ for $t\geq t_r$. \\
\textit{Step 2: the convergence of} $z_i$. Consider a Lyapunov function candidate

\begin{equation*}
V_z= \frac{1}{2}\tilde{z}^T (H\otimes I_3)\tilde{z}
\end{equation*}
which satisfies $V_z\geq 0$ and $V_z=0$ if and only if $\tilde{z}=0$, since $H$ is positive definite according to Lemma \ref{lem_graph}. The time derivative of $V_z$ along (\ref{eq_errobs3}) is computed as

\begin{equation} \label{eq_B0}
\begin{array}{cll}
\dot{V}_z & = & \tilde{z}^T (H\otimes I_3)\dot{\tilde{z}} \\
 {} & = & z_s^T [-\lambda_3\textup{sgn}(z_s - A_0  \tilde{w}) - \textbf{1}_n\otimes{\dot{z}_0}]
\end{array}
\end{equation}
Note that

\begin{equation*}
\begin{array}{cll}
\dot{V}_z & \leq & (\lambda_3 + \gamma_3) \|z_s\|_2 \\
{} & \leq & (\lambda_3 + \gamma_3) \bar{\sigma}(H) \|\tilde{z}\|_2   \\
{} & \leq & c_{w}\sqrt{V_z} \\
{} & {} & \text{with}\, c_w=(\lambda_3 + \gamma_3) \bar{\sigma}(H) \sqrt{\frac{2} {\underline{\sigma}(H)}}
\end{array}
\end{equation*}
where $\underline{\sigma}(H\otimes I_3)= \underline{\sigma}(H)$ and $\bar{\sigma}(H\otimes I_3)= \bar{\sigma}(H)$ are utilized in the above derivations. Invoking the comparison principle \cite{Khalil:02} then implies that $\tilde{z}(t)$ remains bounded for $t<t_r$. For $t \geq t_r$, $A_0\tilde{w}(t)=0$ and~(\ref{eq_B0}) becomes

\begin{equation} \label{eq_B1}
\begin{split}
\dot{V}_z &= z_s^T [-\lambda_3\textup{sgn}(z_s) - \textbf{1}_n\otimes{\dot{z}_0}] \\
&\leq -\lambda_3\|z_s\|_1 + \|\dot{z}_0\|_{\infty}\|z_s\|_1  \\
&\leq -(\lambda_3 - \gamma_3)\|z_s\|_1
\end{split}
\end{equation}
where $\dot{z}_0=\ddot{\omega}_0$ and $\|\ddot{\omega}_0\| \leq \gamma_3$ are used in deriving (\ref{eq_B1}). $\forall x\in\mathbb{R}^n$, Lemma \ref{lem_hardy} can be used to show that
$$\|x\|_2 \leq \|x\|_1 \leq \sqrt{n}\|x\|_2$$
Noting $\underline{\sigma}(H\otimes I_3)= \underline{\sigma}(H)$ and $\bar{\sigma}(H\otimes I_3)= \bar{\sigma}(H)$, one obtains $0.5{\underline{\sigma}(H)}\|\tilde{z}\|_2^2 \leq V_z \leq 0.5{\bar{\sigma}(H)}\|\tilde{z}\|_2^2$ and
\begin{equation*}
\|z_s\|_1 \geq \|z_s\|_2 = \|(H\otimes I_3)\tilde{z}\|_2 \geq \underline{\sigma} (H)\|\tilde{z}\|_2.
\end{equation*}
It then follows from (\ref{eq_B1}) that
\begin{equation*}
\dot{V}_z \leq -c_z \sqrt{V_z}, \; \textup{with} \; c_z= (\lambda_3 - \gamma_3)\underline{\sigma}(H)\sqrt{\frac{2}{\bar{\sigma}(H)}}.
\end{equation*}
which implies, according to the comparison principle , that $\tilde{z}(t)$ (or, equivalently, $z(t)$) is uniformly bounded and $\tilde{z}(t)= 0$ for all $t \geq T_z$, where
\begin{equation*}
T_z= \frac{2\sqrt{V_z(\tilde{z}(0))}}{c_z} + t_r.
\end{equation*}
\textit{Step 3: the convergence of} $v_i$. Similarly, consider the time derivative of $V_v= \tilde{v}^T (H\otimes I_3) \tilde{v}/2$ along (\ref{eq_errobs2}):
\begin{equation} \label{eq_B2}
\dot{V}_v = \tilde{v}^T (H\otimes I_3)\dot{\tilde{v}} =-\lambda_2v_s^T\textup{sgn}^{\beta_2}(v_s) + v_s^T\tilde{z}.
\end{equation}
Recall that (\ref{eq_B1}) implies ${V}_z(\tilde{z}(t)) \leq {V}_z(\tilde{z}(0))$ and thus ${\|\tilde{z}(t)\|_2}\leq{\gamma_z} \triangleq \sqrt{2V_z(\tilde{z}(0))/{\underline{\sigma}(H)}}$. In addition, $\underline{\sigma}(H){\|\tilde{v}\|_2} \leq \|v_s\|_2 \leq \bar{\sigma}(H){\|\tilde{v}\|_2}$ and $0.5{\underline{\sigma}(H)}\|\tilde{v}\|_2^2 \leq V_v \leq 0.5{\bar{\sigma}(H)}\|\tilde{v}\|_2^2$. It then follows from (\ref{eq_B2}) that $\dot{V}_v \leq v_s^T\tilde{z} \leq{\gamma_z \|v_s\|_2} \leq{c_{v0}\sqrt{V_v}}$, where $c_{v0}=\gamma_z\bar{\sigma}(H)\sqrt{{2}/{\underline{\sigma}(H)}}$, and thus
$$V_v(t) \leq \Gamma_v(t) \triangleq \left[0.5c_{v0}t + \sqrt{V_v(\tilde{v}(0))}\right]^2,$$
which implies that $\tilde{v}(t)$ cannot escape in finite time and for $t\in[T_z,\infty)$ (\ref{eq_B2}) reduces to $\dot{V}_v= -\lambda_2v_s^T\textup{sgn}^{\beta_2}(v_s)$. Note that $\forall x\in\mathbb{R}^n$, Lemma \ref{lem_hardy} implies $x^T\textup{sgn}^{\beta_2}(x)= \sum_{i=1}^{n}|x_i|^{1+\beta_2} \geq (\sum_{i=1}^{n}|x_i|^2)^{(1+\beta_2)/2}$. It then follows that
\begin{equation} \label{eq_B3}
\begin{array}{cc}
  \begin{split}
     \dot{V}_v &\leq -\lambda_2{\|v_s\|_2^{1+\beta_2}} \\
     &\leq -c_v V_v^{\frac{1+\beta_2}{2}}
  \end{split}, & c_v= \lambda_2\left[\frac{2\underline{\sigma}^2(H)}{\bar{\sigma}(H)}\right]^{\frac{1+\beta_2}{2}}.
\end{array}
\end{equation}
and hence $\tilde{V}_v(t)\leq \Gamma_v(T_z)$. Noting $\|v\|_2 \leq \|\tilde{v}\|_2 + \sqrt{n}\|v_0\|_2$ and $\|v_0\|_2 \leq \sqrt{3}\gamma_1$, one can deduce that $\|\tilde{v}(t)\|_2,\|v(t)\|_2 \leq \gamma_v$, where $\gamma_v = \sqrt{2\Gamma_v(T_z)/{\underline{\sigma}(H)}} + \sqrt{3n}\gamma_1$. Noting $0 < (1+\beta_2)/2 < 1$ and applying again the comparison principle, it follows that $\tilde{v}(t)\to 0$ in a finite time $T_v \geq T_z$ satisfying $T_v - T_z \leq {2V_v^{(1-\beta_2)/{2}}(\tilde{v}(T_z))}/{c_v(1-\beta_2)}$. Furthermore, an estimation of $T_v$ can be obtained as
\begin{equation*}
T_v \leq T_z + \frac{2\Gamma_v^{(1-\beta_2)/2}(T_z)}{c_v(1-\beta_2)}.
\end{equation*}
\textit{Step 4: the convergence of} $P_i$. The proof can be performed in a manner similar to Step 2. More precisely, we first show the absence of finite escape time for $\tilde{P}(t)$ and then verify the convergence of the reduced system of (\ref{eq_errobs1}) for $t \geq T_v$. To this end, consider a Lyapunov function candidate $V_p= \tilde{P}^T (H\otimes I_4) \tilde{P}/2$ and its time derivative along (\ref{eq_errobs1}):
\begin{equation} \label{eq_B4}
\begin{split}
   \dot{V}_p & = \tilde{P}^T (H\otimes I_4)\dot{\tilde{P}} \\
     & =-\lambda_1 P_s^T\textup{sgn}^{\beta_1}(P_s) + \frac{1}{2} \sum_{i=1}^{n} P_{si}^T[E(\tilde{P}_i)v_i + E(P_0)\tilde{v}_i].
\end{split}
\end{equation}
Since $\tilde{v}(t)$ and $v(t)$ are uniformly bounded, there exists a constant $\gamma_v \geq 0$ such that $\|\tilde{v}(t)\|_2,\|v(t)\|_2 \leq \gamma_v$. Invoking Lemma \ref{lem_E1} and noting $\|P_0\|_2=\|Q_0\|_2=1$, the following inequalities can be shown:
\begin{equation*}
\begin{array}{rll}
P_{si}^TE(\tilde{P}_i)v_i & \leq & \|P_{si}\|_2 \|E(\tilde{P}_i)\|_2 \|v_i\|_2 \\
{} & \leq & \gamma_v \|P_{s}\|_2 \|\tilde{P}_{i}\|_2 \leq \gamma_v \bar{\sigma}(H)\|\tilde{P}\|_2^2, \\
P_{si}^T E(P_0)\tilde{v}_i & \leq & \gamma_v \bar{\sigma}(H)\|\tilde{P}\|_2 \\
{} & \leq & \gamma_v \bar{\sigma}(H)[0.25 + \|\tilde{P}\|_2^2].
\end{array}
\end{equation*}
It then follows from (\ref{eq_B4}) that
$\dot{V}_p \leq c_{p0}V_p + c_{p1}$, where $ c_{p0}={2n\gamma_v \bar{\sigma}(H)}/{\underline{\sigma}(H)}$ and $c_{p1}={n\gamma_v \bar{\sigma}(H)}/{8}$. Invoking the comparison principle leads to $V_p(t)\leq c_{p2}e^{c_{p0}t} - c_{p1}/c_{p0}$, where $c_{p2}=V_p(\tilde{P}(0))+c_{p1}/c_{p0}$, and thus $\tilde{P}(t)$ has no finite escape time. For $t\in [T_v,\infty)$, it follows that $\tilde{v}(t)=0$ and $v_i(t)=v_0$, $i\in\mathbb{I}_n$, and (\ref{eq_B4}) reduces to
\begin{equation*} 
\begin{split}
   \dot{V}_p &= -\lambda_1 P_s^T\textup{sgn}^{\beta_1}(P_s) + \frac{1}{2} \sum_{i=1}^{n} P_{si}^T E(\tilde{P}_i)v_0.
\end{split}
\end{equation*}
In addition, we can deduce by means of (\ref{eq_obsin}) and Lemmas \ref{lem_E1} and \ref{lem_E2} that
\begin{equation*}
\begin{split}
   \sum_{i=1}^{n}P_{si}^T E(\tilde{P}_i)v_0 &= \sum_{i=1}^{n}\sum_{j=1}^{n}a_{ij}(\tilde{P}_i - \tilde{P}_j)^T E(\tilde{P}_i)v_0 \\
   &= -\sum_{i=1}^{n}\sum_{j=1}^{n}a_{ij}\tilde{P}_j^T E(\tilde{P}_i)v_0 =0
\end{split}
\end{equation*}
As a result, derivations similar to (\ref{eq_B3}) can be used to show that $\dot{V}_p \leq -c_p V_p^{(1+\beta_1)/{2}}$ holds for $t\in [T_v,\infty)$, where $c_p= \lambda_1 \left[{2\underline{\sigma}^2(H)}/{\bar{\sigma}(H)}\right] ^{(1+\beta_1)/{2}}$. Therefore, $\tilde{P}(t)$ (and thus $P(t)$) is uniformly bounded and $\tilde{P}(t)\to 0$ in a finite time $T_p \geq T_v$ satisfying  $T_p - T_v \leq {2V_p^{(1-\beta_1)/{2}}(\tilde{P}(T_v))}/{c_p(1-\beta_1)}$. An estimation of $T_p$ is derived as
\begin{equation} \label{eq_Tp}
T_p \leq T_v + \frac{2\left[c_{p2}e^{c_{p0}T_v} - c_{p1}/c_{p0}\right]^\frac{1-\beta_1}{2}}{c_p(1-\beta_1)}.
\end{equation}
Summarizing all the above arguments follows the uniform boundedness of $(P_i(t),v_i(t),z_i(t))$ and $(P_i(t), v_i(t),\\ z_i(t))=(Q_{0}(t),\omega_{0}(t),\dot{\omega}_{0}(t))$ for $t \geq T_p$, $i\in{\mathbb I_n}$, thus concluding the results of Theorem \ref{thm_obs}.

\section{Proof of Lemma \ref{lem_kappa}}
\label{app3}       
Given $Q=[\eta,q^T]^T \in\mathbb{R}^4$, it is clear that $\bar{\kappa}_1(Q,\alpha)$ is continuous with respect to $Q$ if $\eta \neq \|Q\|_2$. In addition, the following identities are straightforward
\begin{equation*}
\begin{array}{c}
\eta^2 + \|q\|_2^2 =\|Q\|_2^2, \\
2\|Q\|_2(\|Q\|_2 - \eta)=(\|Q\|_2-\eta)^2 + \|q\|_2^2.
\end{array}
\end{equation*}
It then follows that $\sqrt{2\|Q\|_2(\|Q\|_2 -\eta)} \geq \|q\|_2$ and thus

\begin{equation*}
\| \bar{\kappa}_1(Q,\alpha) \|_2 \leq {\|q\|_2^{1-\alpha}}= \left(\sqrt{\|Q\|_2^2 - \eta^2}\right)^{1-\alpha}
\end{equation*}
Therefore, we have

\begin{equation*}
\| \bar{\kappa}_1(Q,\alpha) \|_2 \leq \|Q\|_2^{1-\alpha} \,\text{and}\,\lim_{\eta\to \|Q\|_2}\bar{\kappa}_1(Q,\alpha)=0
\end{equation*}
which implies that $\bar{\kappa}_1(Q,\alpha)$ is continuous at $\eta= \|Q\|_2$ and thus on $\mathbb{R}^4$.

\section{Proof of Theorem \ref{thm_GFTC2}}
\label{app4}

Since $u_i(\hat{\mathbf{x}}_{fi}(t))= u_i(\mathbf{x}_{fi}(t))$ for $t \geq T_p$, we only need to show that there is no finite escape time for the closed-loop trajectory $(Q_{i0}(t),\omega_{i0}(t))$. Note that $Q_{i0}(t)\in\mathbb{S}^3$ is trivially bounded. Next, consider the positive-definite function $V_i= \omega_{i0}^T J_i \omega_{i0}/2$. For $\hat{\mathbf{x}}_{fi}\in \hat{C}_{fi}$, its time derivative along (\ref{eq_errdyn}) and (\ref{eq_GFTC2}) is computed as
\begin{equation} \label{eq_D1}
\dot{V}_i= \omega_{i0}^T[\hat{u}_{fi} - u_{fi} -k_{pi}\bar{\kappa}_1(h_i \hat{Q}_{i0},1-\alpha_{pi}) - k_{di}\textup{sat}_{\alpha_{di}}(\hat{\omega}_{i0})].
\end{equation}
Noting that $\|R(\hat{Q}_{i0})\|_2 = \|\hat{Q}_{i0}\|_2^2$ and $\|\hat{Q}_{i0}\|_2= \|P_i\|_2$, Assumption \ref{assum_leader} and Lemma \ref{lem_kappa} together with the uniform boundedness of $(P_i,v_i,z_i)$ can then be used to show the uniform boundedness of all terms involved in the bracket of (\ref{eq_D1}). Hence, there exists a constant $\gamma_u>0$ such that $\dot{V}_i \leq \gamma_u \|\omega_{i0}\|_2 \leq c_i\sqrt{V_i}$, where $c_i=\gamma_u\sqrt{2/\underline{\sigma}(J_i)}$. Invoking the comparison principle and recalling the continuity of $\omega_{i0}(t)$ over jumps, it follows that $\omega_{i0}(t)$ is bounded for $0\leq t<T_p$. When $t \geq T_p$, $u_i(\hat{\mathbf{x}}_{fi}(t))= u_i(\mathbf{x}_{fi}(t))$ and the closed-loop equations reduce to a globally finite-time stable system, according to Theorem 3.1 in \cite{Gui:16SCL}. Therefore, the statements of Theorem \ref{thm_GFTC2} follow.

\end{document}